\documentclass[12pt]{amsart}
\usepackage{amsfonts,amssymb,amsmath,amscd,amsthm, amstext,bm,color}
\usepackage[colorlinks, linkcolor=blue, citecolor=blue, urlcolor=blue, pagebackref,hypertexnames=false]{hyperref}
\usepackage{enumerate}
\usepackage{fancyhdr}
\usepackage{geometry}
\usepackage{graphicx}
\usepackage{epstopdf}
\usepackage{mathrsfs}
\usepackage{txfonts}
\usepackage{tikz}
\usetikzlibrary{arrows.meta, positioning, decorations.pathreplacing}
\usepackage{url}
\usepackage{comment}
\allowdisplaybreaks
\numberwithin{equation}{section}
\newtheorem{theorem}{Theorem}[section]
\newtheorem{proposition}[theorem]{Proposition}
\newtheorem{lemma}[theorem]{Lemma}
\newtheorem{corollary}[theorem]{Corollary}
\newtheorem{remark}[theorem]{Remark}
\newcommand{\R}{\mathbb R}
\newcommand{\Xnorm}[1]{  \|#1   \|_{\mathbb L^2(\R_+;\mathrm{BMO}(\mu_\kappa))}}
\newcommand{\Bnorm}[1]{\left\|#1\right\|_{\mathrm{BMO}_2(\mu_\kappa;\mathcal H)}}

\title[Global maximal regularity for the Dunkl heat equation]
{Global-in-time maximal regularity for the Dunkl heat equation in BMO}
\author{Xuejing Huo}
\address{Xuejing Huo, School of Mathematical and Physical Sciences,
Macquarie University, Sydney 2109, NSW, Australia}
\email{xuejing.huo@hdr.mq.edu.au}

\author{Ming-Yi Lee}
\address{Department of Mathematics, National Central University,
Chung-Li 320, Taiwan, Republic of China}
\email{mylee@math.ncu.edu.tw}

\author{Chaojie Wen}
\address{Department of Mathematics, Sun Yat-sen University, Guangzhou, 510275, China}
\email{wenchj@mail2.sysu.edu.cn}

\author{Liangchuan Wu}
\address{Liangchuan Wu, School of Mathematical Sciences, Anhui University, Hefei 230601, P.R.~China}
\email{wuliangchuan@ahu.edu.cn}
\subjclass[2020]{Primary 42B35; Secondary 35K15, 42B20}
\keywords{Dunkl heat equation, maximal regularity, bounded mean oscillation, square functions}

\begin{document}

\begin{abstract}
We establish global-in-time maximal regularity for the Cauchy problem of the Dunkl heat equation in BMO, without assuming reflection invariance of the external force.
For a Hilbert-valued BMO function $F$, we integrate the spatial mean-square oscillation
of the Dunkl gradient of its heat extension over all $t>0$ on each fixed Euclidean ball.
The supremum of these integrals is comparable to the squared BMO seminorm of $F$.
In particular, the upper estimates imply time-weighted energy estimates for the zero-initial-value solution and maximal regularity for nonzero initial values whose weak Dunkl derivatives belong to BMO.
\end{abstract}
\maketitle

\section{Introduction and main results}\label{sec:1}

For the classical heat equation
\begin{equation}\label{e:cl}
 \partial_tu-\Delta u=f\quad\text{in }\R^N\times(0,\infty),\qquad u(\cdot,0)=0,
\end{equation}
when $f\in L^p(\R_+;L^q(\R^N))$ with $1<p,q<\infty$, the solution
satisfies the maximal regularity estimate
\begin{equation}\label{e:lp}
 \|\partial_tu\|_{L^p(\R_+;L^q)}  +\|\Delta u\|_{L^p(\R_+;L^q)}   \leq C\|f\|_{L^p(\R_+;L^q)}.
\end{equation}
Here the $L^p$ norm is taken in time and the $L^q$ norm in space.
Such estimates provide the linear bounds needed for fixed-point arguments
in nonlinear parabolic problems.

When $q=2$, the spatial space $L^2(\R^N)$ is a Hilbert space.
Since the heat semigroup is bounded and analytic on this space,
de Simon's theorem \cite{dS} implies \eqref{e:lp}
for every $1<p<\infty$;
see also \cite[Corollary~1.7]{KW04}.
For general $1<q<\infty$, the spatial space $L^q(\R^N)$
has the unconditional martingale differences (UMD) property.
Criteria for maximal regularity on UMD spaces have been developed
through operator-valued singular integrals and functional calculus;
see \cite{W,DHP,KW04,HNVW}.
For applications to nonlinear parabolic problems, we refer to \cite{PS}.

At the spatial endpoint $q=\infty$, \eqref{e:lp} fails; see \cite[Example~17.4.3]{HNVW}.
The space $\mathrm{BMO}(\R^N)$ is a natural substitute for $L^\infty(\R^N)$ at this endpoint.
For maximal regularity, one must also specify how
the spatial mean oscillation is integrated in time.
Ogawa and Shimizu \cite{OS} established maximal regularity for \eqref{e:cl}
with a seminorm equivalent to
$$
 \sup_{B=B(c,r_B)}
 \bigg (\frac1{|B|}\int_0^{r_B^2}\int_B
\bigg |f(x,t)-\frac1{|B|}\int_B f(y,t)\,dy\,\bigg|^2\,dx\,dt\bigg)^{1/2}.
$$
In the global result of \cite{OS}, the solution exists for all $t>0$.
For each ball $B$, the time integral in the seminorm runs only over $(0,r_B^2)$
(see Definition~1.1, (1.10) and Corollary~1.3 therein).
Parabolic tent spaces likewise couple the time and spatial scales;
see \cite{AF,AH,AKMP,AMP} for their maximal regularity theory.

Recently, Duong, Li, the fourth author and Yan \cite{DLWY} established maximal regularity
for the classical heat equation with the time integral taken over $(0,\infty)$
in the BMO seminorm above.
 For $t>r_B^2$, the Duhamel integral over $0<s<t$ splits into the ranges
$0<t-s<r_B^2$ and $t-s\geq r_B^2$. In the latter range,
spatial regularity of the heat kernel yields additional decay in $t-s$.
Analogous estimates were also established for Schr\"odinger
operators with reverse H\"older potentials in the corresponding
adapted BMO spaces.

In this paper, we investigate global-in-time maximal regularity for the Dunkl
heat equation in BMO.
Let $R$ be a finite reduced root system in $\R^N$, let $R_+$ be a choice of positive roots,
and let $G$ be generated by the reflections
$$
 \sigma_\alpha x=x-2\frac{\langle\alpha,x\rangle}{\|\alpha\|^2}\alpha,
 \qquad \alpha\in R .
$$
Here $\langle\cdot,\cdot\rangle$ and $\|\cdot\|$ denote
the Euclidean inner product and norm on $\R^N$, respectively.
For a nonnegative $G$-invariant multiplicity function $\kappa$, write
$$
 d\mu_\kappa(x)=\prod_{\alpha\in R_+}|\langle\alpha,x\rangle|^{2\kappa(\alpha)}\,dx.
$$
The associated Dunkl Laplacian is
\begin{equation}\label{e:lap}
 \Delta_\kappa\varphi(x)=\Delta\varphi(x)+2\sum_{\alpha\in R_+}\kappa(\alpha)
 \left(
 \frac{\langle\alpha,\nabla\varphi(x)\rangle}{\langle\alpha,x\rangle}
 -\frac{\|\alpha\|^2}{2}\frac{\varphi(x)-\varphi(\sigma_\alpha x)}{\langle\alpha,x\rangle^2}
 \right).
\end{equation}
The operator is symmetric on $C_c^\infty(\R^N)$ with respect to $\mu_\kappa$.
We use the nonnegative self-adjoint extension of $-\Delta_\kappa$ associated with
the closed quadratic form in \eqref{e:form}, and write $h_t(x,y)$ for its heat kernel.
For the Dunkl heat equation, Deleaval and Kriegler \cite{DK}
established maximal regularity results via H\"ormander spectral
multipliers on UMD lattice-valued $L^p(\mu_\kappa)$, $1<p<\infty$.

We study maximal regularity at the spatial BMO endpoint
for the Cauchy problem
\begin{equation}\label{e:pde}
 \partial_tu-\Delta_\kappa u=f
 \quad\text{in }\R^N\times(0,\infty),\qquad u(\cdot,0)=0.
\end{equation}
For the Duhamel solution
\begin{equation}\label{e:duh}
 u(x,t)=\int_0^t e^{(t-s)\Delta_\kappa}f(\cdot,s)(x)\,ds,
\end{equation}
we seek to estimate $\partial_tu$ and $\Delta_\kappa u$ in terms of $f$
using their spatial mean oscillations with respect to $\mu_\kappa$,
integrated in time over $(0,\infty)$ for each fixed Euclidean ball.
In \cite{OS} the corresponding time interval is $(0,r_B^2)$ for $B=B(c,r_B)$.

For a spatial function $F$ and a Euclidean ball $B\subset\R^N$, write
$$
 F_B=\frac1{\mu_\kappa(B)}\int_B F(x)\,d\mu_\kappa(x).
$$
At time $t$, the spatial average of $f$ over $B$ with respect to $\mu_\kappa$ is
$$
 (f(\cdot,t))_B
 =\frac1{\mu_\kappa(B)}\int_B f(x,t)\,d\mu_\kappa(x).
$$
Following \cite{DLWY}, define
\begin{equation}\label{e:bmo}
 \Xnorm{f}^2
 =\sup_B\frac1{\mu_\kappa(B)}
 \int_0^\infty\int_B  \left|f(x,t)-(f(\cdot,t))_B\right|^2\,d\mu_\kappa(x)\,dt,
\end{equation}
where the supremum is over Euclidean balls.
We write $\mathbb L^2(\R_+;\mathrm{BMO}(\mu_\kappa))$
for the space of functions with finite seminorm \eqref{e:bmo}.
In this seminorm, we integrate in time before taking the supremum over balls.
For the usual Bochner seminorm on
$L^2(\R_+;\mathrm{BMO}(\mu_\kappa))$,
the supremum is taken at each time before integration.
 The seminorm \eqref{e:bmo}  does not change when a function of time alone is added to $f$.

We say that $u$ is a weak solution of
$\partial_tu-\Delta_\kappa u=f$ with initial value $u_0$ if
\begin{equation}\label{e:weak}
 \int_0^\infty\int_{\R^N}
 u(-\partial_t\varphi-\Delta_\kappa\varphi)\,d\mu_\kappa\,dt
 = \int_0^\infty\int_{\R^N}f\varphi\,d\mu_\kappa\,dt
 +\int_{\R^N}u_0(x)\varphi(x,0)\,d\mu_\kappa(x)
\end{equation}
for every $\varphi\in C_c^\infty(\R^N\times[0,\infty))$.
The condition $\Delta_\kappa u\in\mathbb L^2(\R_+;\mathrm{BMO}(\mu_\kappa))$
means that this distribution is represented by a locally integrable function
on $\R^N\times(0,\infty)$ with finite seminorm \eqref{e:bmo}.

\begin{theorem}\label{t:main}
Let $f\in L^1_{\mathrm{loc}}(\R^N\times[0,\infty),d\mu_\kappa\,dt)$
with $\Xnorm{f}<\infty$.
The function $u$ defined by \eqref{e:duh} is the unique weak solution
of \eqref{e:pde} satisfying
\begin{equation}\label{e:class}
 u\in C([0,\infty);L^2_{\mathrm{loc}}(\mu_\kappa)),\qquad
 \Delta_\kappa u\in\mathbb L^2(\R_+;\mathrm{BMO}(\mu_\kappa)).
\end{equation}
Moreover,
\begin{equation}\label{e:mr}
 \Xnorm{\partial_tu}+\Xnorm{\Delta_\kappa u}\leq C\Xnorm{f},
\end{equation}
where $C>0$ is a constant independent of $f$.
\end{theorem}

Neither the external force nor the solution is required to be $G$-invariant.
In the proof, the time variable is absorbed into the Hilbert space $L^2(\R_+)$.
For a separable Hilbert space $\mathcal H$, let
$$
 \Bnorm{F}
 =\sup_B\left(\frac1{\mu_\kappa(B)}
 \int_B  \left \|F(x)-F_B\right\|_{\mathcal H}^2\,d\mu_\kappa(x)\right)^{1/2}
 \qquad  \text{for}\  \   F\in L^2_{\mathrm{loc}}(\mu_\kappa;\mathcal H).
$$
In the scalar case, we omit $\mathcal H$ and write $\mathrm{BMO}(\mu_\kappa)=\mathrm{BMO}_2(\mu_\kappa)$.

By Lemma~\ref{l:norm}, the difference
$f(x,t)-(f(\cdot,t))_{B(0,1)}$, regarded as an $L^2(\R_+)$-valued function of $x$,
belongs to $\mathrm{BMO}_2(\mu_\kappa;L^2(\R_+))$ with seminorm $\Xnorm{f}$.
Its $L^2(\R_+)$ norm is locally square-integrable in $x$. Thus the local integrability
assumption in Theorem~\ref{t:main} amounts to $(f(\cdot,t))_{B(0,1)}\in L^1(0,T)$ for every $T>0$.
Adding $c(t)$ to $f$, with $c\in L^1_{\mathrm{loc}}([0,\infty))$,
adds $\int_0^t c(s)\,ds$ to $u$ and preserves \eqref{e:mr}.

For strongly measurable BMO-valued $f$, the usual Bochner seminorm satisfies
$$
 \|f\|_{L^2(\R_+;\mathrm{BMO}(\mu_\kappa))}^2
 =\int_0^\infty\sup_B\frac1{\mu_\kappa(B)}
 \int_B|f(x,t)-(f(\cdot,t))_B|^2\,d\mu_\kappa(x)\,dt.
$$
Taking the supremum outside the time integral cannot increase this quantity:
$$
 \Xnorm{f}\leq\|f\|_{L^2(\R_+;\mathrm{BMO}(\mu_\kappa))}.
$$
Estimate \eqref{e:mr} fails with this stronger seminorm on both sides.
The case $\kappa=0$ was noted in \cite[Section~1]{DLWY}.
Proposition~\ref{p:bm} proves the failure on every finite time interval
for smooth, compactly supported external forces.

\begin{remark}\label{r:orb}
The image of $B=B(c,r)$ under $\sigma\in G$ is the Euclidean ball
$\sigma B=\{\sigma x:x\in B\}=B(\sigma c,r)$.
The estimate \eqref{e:mr} fails in general if $B$ is replaced by
the orbit ball $\bigcup_{\sigma\in G}\sigma B$ in all three seminorms.
In rank one, for each $\kappa\geq0$, the odd external forces
supported in $\R\times(0,1)$ constructed in Proposition~\ref{p:orb}
are counterexamples.
\end{remark}

 The reflection terms in \eqref{e:lap} explain
why reflected balls enter the proof of Theorem~\ref{t:main}.
Let $x$ lie outside the reflecting hyperplanes, and let $\varphi$ be a smooth function vanishing near $x$.
Then \eqref{e:lap} reduces to
$$
 \Delta_\kappa\varphi(x)
 =\sum_{\alpha\in R_+}
 \frac{\kappa(\alpha)\|\alpha\|^2}{\langle\alpha,x\rangle^2}
 \varphi(\sigma_\alpha x).
$$
The reflection terms thus couple $x$ to its reflected points,
whereas \eqref{e:bmo} measures oscillation on individual Euclidean balls.

For all $x,y\in\R^N$ and $t>0$, the Dunkl heat-kernel bounds in
\cite[Theorem~3.1]{DH1} contain the Gaussian factor $e^{-c_0d(x,y)^2/t}$, where $c_0>0$ and
$$
 d(x,y)=\min_{\sigma\in G}\|x-\sigma y\|.
$$
If $x\in B=B(c,r)$ and $y\in\sigma B$, then
$$
 d(x,y)\leq\|x-\sigma^{-1}y\|<2r,
$$
regardless of $\|c-\sigma c\|$.
Consequently, the Gaussian factor need not be small
near a distant reflected ball.
To estimate the contribution from $\sigma B$, we write
$$
 f(y,t)-(f(\cdot,t))_B
 =f(y,t)-(f(\cdot,t))_{\sigma B}
 +(f(\cdot,t))_{\sigma B}-(f(\cdot,t))_B.
$$
The first term is estimated by the mean oscillation on $\sigma B$
in \eqref{e:bmo}.
By repeatedly doubling the radii of the two balls, we obtain
$$
 \bigl\|(f(\cdot,t))_{\sigma B}-(f(\cdot,t))_B
 \bigr\|_{L^2(\R_+,dt)}
 \leq C\log\left(2+\frac{\|c-\sigma c\|}{r}\right)\Xnorm{f}.
$$
This logarithmic growth occurs in the rank-one example
of Proposition~\ref{p:ex}.

To control the logarithm above, we use the factor
$(1+\|x-y\|/\sqrt\tau)^{-2}$ in the bound for $h_\tau(x,y)$ at time $\tau>0$
(see \cite[Theorem~3.1]{DH1}). In the Duhamel formula \eqref{e:duh}, $\tau=t-s$.
If $x\in B$ and $y\in\sigma B$, we have
$$
 \|c-\sigma c\|
 \leq \|c-x\|+\|x-y\|+\|y-\sigma c\|
 \leq 2r+\|x-y\|.
$$
Consequently,
\begin{align*}
 \left(1+\frac{\|x-y\|}{\sqrt\tau}\right)^{-2}
 \log\left(2+\frac{\|c-\sigma c\|}{r}\right)
 &\leq
 \frac{\log(4+\|x-y\|/r)}
      {(1+\|x-y\|/\sqrt\tau)^2}\\
 &\leq
 \frac{\log4+|\log(\sqrt\tau/r)|
       +\log(1+\|x-y\|/\sqrt\tau)}
      {(1+\|x-y\|/\sqrt\tau)^2}\\
 &\leq C\left(1+\left|\log\frac{\sqrt\tau}{r}\right|\right),
\end{align*}
with a bound independent of $\|c-\sigma c\|$.
Lemma~\ref{l:mean} uses this inequality to estimate
the spatial integrals involving the heat kernel.
Euclidean decay is also used to control differences of reflected averages
in \cite[Proposition~4.5]{DH4}.

The Euclidean and orbit distances also enter the kernel estimates
in the Dunkl singular-integral theory developed by
Tan, Han, Han, the second author and Li \cite{THHLL,THp}.
By $H^1$--BMO duality, Theorems~1.27 and 4.6 of \cite{THp} imply that scalar Dunkl Calder\'on--Zygmund operators
which vanish on constant functions are bounded on
$\mathrm{BMO}(\mu_\kappa)$.
To estimate $\Delta_\kappa u$, we apply Proposition~\ref{p:cz} with Hilbert space
$L^2(\R_+)$ to the truncations \eqref{e:cut}. Their spatial BMO seminorms
are bounded uniformly in the truncation endpoints. In Section~\ref{sec:3},
estimates for the spatial means complete the construction of the weak solution
and the proof of uniqueness in \eqref{e:class}.

 We next   compare the BMO seminorm of
a spatial function $F$ with the mean oscillation of the Dunkl gradient
of its heat extension $e^{t\Delta_\kappa}F$.

Kaiser and Weis \cite[Theorem~4.2 and Remark~4.4]{KW}
obtained equivalent BMO norms using vector-valued wavelet transforms
and reconstruction formulas.
Dziuba\'nski and Hejna \cite{DH2} proved upper and lower $L^p$ bounds
for Dunkl square functions using Hilbert-valued kernels
and reproducing formulas.
 For scalar functions in orbit BMO, Almeida et al.~\cite{ABFR}
studied square functions defined by time derivatives of the Dunkl heat semigroup,
assuming that the square functions are finite almost everywhere.

 Dziuba\'nski and Hejna \cite[Corollary~4.15]{DH4} proved a Carleson estimate
for Dunkl convolutions of BMO functions with smooth, compactly supported kernels of integral zero.

 To state the corresponding estimate for the Poisson extension, let
$\nabla=(\partial_1,\ldots,\partial_N)$ denote the Euclidean gradient in $x$
and $\nabla_\kappa=(T_1,\ldots,T_N)$ the Dunkl gradient, where
$$
 T_jF(x)=\partial_jF(x)+\sum_{\alpha\in R_+}\kappa(\alpha)\alpha_j
 \frac{F(x)-F(\sigma_\alpha x)}{\langle\alpha,x\rangle},
 \qquad 1\leq j\leq N.
$$
For a scalar BMO function $F$ and its Dunkl Poisson extension
 $u(\cdot,t)=e^{-t\sqrt{-\Delta_\kappa}}F$,
\cite[Corollary~4.17]{DH4} and \cite[Proposition~7.4]{JL} imply
$$
 \|F\|_{\mathrm{BMO}(\mu_\kappa)}^2
 \simeq \sup_{B=B(c,r)}\frac1{\mu_\kappa(B)}
 \int_0^r\int_B t\, \Big  (|\partial_tu|^2+\sum_{j=1}^N|T_ju|^2\Big)
 \,d\mu_\kappa\,dt.
$$
Dziuba\'nski and \L{}y\.zwa \cite[Theorem~1.1]{DL} proved that every smooth Dunkl
harmonic function in the upper half-space satisfying this Carleson condition
is the Poisson extension of a BMO function.

In the Poisson estimate above, $t$ ranges from $0$ to $r$ for each ball $B=B(c,r)$.
For the heat extension $e^{t\Delta_\kappa}F$,
we subtract the spatial mean of the Dunkl gradient on $B$ at each $t$
and integrate its squared oscillation over $(0,\infty)$.
Again, let $\mathcal H$ be a separable Hilbert space and let
$F\in\mathrm{BMO}_2(\mu_\kappa;\mathcal H)$. By Lemma~\ref{l:conv},
the heat-kernel integral defining $e^{t\Delta_\kappa}F$ converges for every $t>0$.
The Dunkl gradient of this extension takes values in $\mathcal H^N$, the product of $N$
copies of $\mathcal H$, with
$\|(v_1,\ldots,v_N)\|_{\mathcal H^N}^2=\sum_{j=1}^N\|v_j\|_{\mathcal H}^2$.
With this norm, \eqref{e:bmo} reads
$$
 \Xnorm{\nabla_\kappa e^{t\Delta_\kappa}F}^2
 =\sup_B\frac1{\mu_\kappa(B)}\int_0^\infty\int_B
 \bigl\|\nabla_\kappa e^{t\Delta_\kappa}F(x)
       -(\nabla_\kappa e^{t\Delta_\kappa}F)_B\bigr\|_{\mathcal H^N}^2
 \,d\mu_\kappa(x)\,dt.
$$
Without this subtraction, the integral diverges even for
$F(x)=\operatorname{sgn}x$ on $\R$ with $\kappa=0$
(Remark~\ref{r:ctr}).

\begin{theorem}\label{t:sq}
For every $F\in\mathrm{BMO}_2(\mu_\kappa;\mathcal H)$,
\begin{equation}\label{e:dsq}
 C^{-1}\Bnorm{F}\leq\Xnorm{\nabla_\kappa e^{t\Delta_\kappa}F}\leq C\Bnorm{F}.
\end{equation}
Moreover,
\begin{equation}\label{e:sq}
 \Xnorm{\nabla e^{t\Delta_\kappa}F}\leq C\Bnorm{F}.
\end{equation}
For each $\alpha\in R_+$ with $\kappa(\alpha)>0$,
\begin{equation}\label{e:rsq}
 \biggl\|\frac{e^{t\Delta_\kappa}F-(e^{t\Delta_\kappa}F)\circ\sigma_\alpha}
 {\langle\alpha,\cdot\rangle}\biggr\|_{\mathbb L^2(\R_+;\mathrm{BMO}(\mu_\kappa))}\leq C\Bnorm{F}.
\end{equation}
The constant is independent of $\mathcal H$.
\end{theorem}

 For the lower bound in \eqref{e:dsq}, the identity \eqref{e:rep}
allows us to estimate the oscillation of
$e^{2a\Delta_\kappa}F-e^{2b\Delta_\kappa}F$ by that of the Dunkl heat
gradient on $a<t<b$. On each fixed ball, the first heat extension
converges to $F$ in $L^2$ as $a\to0+$, and the mean oscillation of
the second tends to zero as $b\to\infty$.
Appendix~\ref{sec:a} contains a second proof of the upper bounds and of
the truncation estimate in Proposition~\ref{p:mr}, following \cite{DLWY}.
For the part of the function outside $4B$, where $B=B(c,r)$, we use Euclidean decay for $0<\tau<r^2$
and use the spatial derivative to compare the kernel at $x\in B$ and at $c$
for $\tau\geq r^2$.
The short-time bound \eqref{e:short} holds without subtracting the spatial
mean of the Dunkl heat gradient.

 The Euclidean gradient and the reflection quotients in
\eqref{e:sq} and \eqref{e:rsq} are the two terms in the energy form
\eqref{e:form}. We apply these bounds to the normalized external force
from Lemma~\ref{l:norm}, with $\mathcal H=L^2(\R_+,ds)$.
The Cauchy--Schwarz inequality in the Duhamel integral over $0<s<t$
then yields the estimates in Corollary~\ref{c:en} for $t^{-1/2}\nabla u$
and the reflection quotients of $t^{-1/2}u$.

 For nonzero initial data we must also estimate $e^{t\Delta_\kappa}u_0$.
If $u_0\in L^2_{\mathrm{loc}}(\mu_\kappa)$ and its weak Dunkl derivatives
belong to $\mathrm{BMO}_2(\mu_\kappa)$, then $u_0$ has polynomial
mean-square growth by Lemma~\ref{l:grow}. The heat-kernel integral
therefore converges for $t>0$. By Lemma~\ref{l:init} and the identity
$\Delta_\kappa=\sum_{j=1}^N T_j^2$,
$$
 \partial_t e^{t\Delta_\kappa}u_0
 =\Delta_\kappa e^{t\Delta_\kappa}u_0
 =\sum_{j=1}^N T_j e^{t\Delta_\kappa}(T_j u_0).
$$
We apply the upper bound in \eqref{e:dsq} to each $T_j u_0$
and Theorem~\ref{t:main} to the Duhamel term.

\begin{corollary}\label{c:init}
Let $f$ be as in Theorem~\ref{t:main} and let
$u_0\in L^2_{\mathrm{loc}}(\mu_\kappa)$ have weak Dunkl derivatives
$T_j u_0\in\mathrm{BMO}_2(\mu_\kappa)$, $j=1,\ldots,N$.
The function
$$
 u(x,t)=e^{t\Delta_\kappa}u_0(x)
 +\int_0^t e^{(t-s)\Delta_\kappa}f(\cdot,s)(x)\,ds
$$
is the unique weak solution of $\partial_tu-\Delta_\kappa u=f$
satisfying \eqref{e:class} and $u(\cdot,0)=u_0$.
Moreover,
\begin{equation}\label{e:imr}
 \Xnorm{\partial_tu}+\Xnorm{\Delta_\kappa u}
 \leq C\bigg(\Xnorm{f}
 +\sum_{j=1}^N\|T_j u_0\|_{\mathrm{BMO}_2(\mu_\kappa)}\bigg),
\end{equation}
where $C>0$ is independent of $f$ and $u_0$.
\end{corollary}

 When $\kappa=0$, the assumption on the initial value becomes
$\nabla u_0\in\mathrm{BMO}(\R^N)$. Ogawa and Shimizu
\cite[Corollary~1.3]{OS} use the time interval $(0,r_B^2)$ on each ball
$B=B(c,r_B)$. Taking $f=0$ in \eqref{e:imr}, we obtain the estimate for
$\partial_t e^{t\Delta}u_0$ over $(0,\infty)$ on each ball
suggested in \cite[Remark~2.1]{DLWY}.

 In rank one with $\kappa>0$, an external force supported near $-a$ can
produce a response of order $t^2$ at $a>0$ for the Dunkl heat equation,
while the classical heat response to the same force is exponentially small there.
Take $f(x,t)=\mathbf1_{(0,1)}(t)\varphi(x)$, where $\varphi$ is smooth,
nonnegative, supported in $(-5a/4,-3a/4)$, and satisfies $\varphi(-a)=1$.
For zero initial value, Proposition~\ref{p:resp} shows that, as $t\to0^+$,
$$
 u_\kappa(a,t)=\frac{\kappa}{2a^2}t^2+O(t^3),
 \qquad
 0\leq u_{\mathrm{cl}}(a,t)
 \leq C\|\varphi\|_\infty\,t\exp(-ca^2/t).
$$
The classical Gaussian decays exponentially since $|a-y|\geq7a/4$ on the support of $\varphi$.
Although $\varphi$ vanishes near $a$, the reflection term in \eqref{e:lap} yields
$\Delta_\kappa\varphi(a)=\kappa\varphi(-a)/a^2=\kappa/a^2$, which accounts for the $t^2$ term.

 The same reflection term produces an inverse-square potential on the half-line.
Subtracting the equations at $r$ and $-r$, we find that
$w(r,t)=u(r,t)-u(-r,t)$ satisfies the Bessel heat equation
$$
 \partial_tw-\partial_r^2w-\frac{2\kappa}{r}\partial_rw
 +\frac{2\kappa}{r^2}w
 =f(r,t)-f(-r,t),\qquad r>0.
$$
Extend the right-hand side to an odd function on $\R$. Its Duhamel solution is odd,
so its derivative and its quotient by $x$ are even. Theorem~\ref{t:main}
and Corollary~\ref{c:en} then yield maximal regularity and the bounds for
the even extensions of $t^{-1/2}\partial_rw$ and $t^{-1/2}w/r$
stated in Corollary~\ref{c:odd}. For maximal $L^p$-regularity of Bessel operators, see \cite{BL}.

This paper is organized as follows. Section~\ref{sec:2} contains the heat-kernel and BMO
estimates and the counterexamples for the Bochner and orbit-ball seminorms.
We prove Theorem~\ref{t:main} in Section~\ref{sec:3} and the heat-gradient
and time-weighted energy estimates in Section~\ref{sec:4}.
Section~\ref{sec:5} treats nonzero initial data. Section~\ref{sec:6} studies
the response at a reflected point and the half-line equation.
Appendix~\ref{sec:a} contains the direct proofs of Proposition~\ref{p:mr}
and the upper bounds in Theorem~\ref{t:sq}.

\section{Heat kernels and BMO estimates}\label{sec:2}

In this section, we prove the Hilbert-valued BMO estimate used for maximal regularity
and for the heat gradient (Proposition~\ref{p:cz}). We compare the averages
on $B$ and $\sigma B$ and use decay in $\|x-y\|$ to control their difference
in the kernel integral.

\subsection{Heat kernels}\label{s:heat}

Recall that the Dunkl Laplacian is $\sum_{j=1}^NT_j^2$.
We say that $u\in L^2_{\mathrm{loc}}(\mu_\kappa)$ has weak derivative $T_ju=g_j$, if
$$
 \int_{\R^N}g_j\varphi\,d\mu_\kappa
 =-\int_{\R^N}uT_j\varphi\,d\mu_\kappa,
 \qquad  \text{for each}\   \  \varphi\in C_c^\infty(\R^N).
$$
The nonnegative self-adjoint operator $-\Delta_\kappa$ is associated with the closure of the form
\begin{equation}\label{e:form}
 \langle-\Delta_\kappa v,v\rangle
 =\int_{\R^N}\|\nabla v\|^2\,d\mu_\kappa
 +\sum_{\alpha\in R_+}\frac{\kappa(\alpha)\|\alpha\|^2}{2}
 \int_{\R^N}\left|\frac{v-v\circ\sigma_\alpha}{\langle\alpha,x\rangle}\right|^2\,d\mu_\kappa.
\end{equation}
We refer to \cite{DH1,DH2} for more details about the Dunkl Laplacian and its heat semigroup.

The ball volumes satisfy
\begin{equation}\label{e:vol}
 \mu_\kappa(B(x,r))
 \simeq r^N\prod_{\alpha\in R_+}
 \bigl(|\langle\alpha,x\rangle|+r\|\alpha\|\bigr)^{2\kappa(\alpha)}.
\end{equation}
Consequently, for $\lambda\geq1$,
$$
 \lambda^N\mu_\kappa(B(x,r))
 \lesssim\mu_\kappa(B(x,\lambda r))
 \lesssim\lambda^{\mathbf N}\mu_\kappa(B(x,r)),
$$
where  $\mathbf N=N+2\sum_{\alpha\in R_+}\kappa(\alpha)$ is the homogeneous dimension.

The heat kernel $h_t(x,y)$ is positive, smooth and symmetric in $x,y$, with
$$
 \int_{\R^N}h_t(x,y)\,d\mu_\kappa(y)=1,\qquad t>0.
$$
After using symmetry and comparing ball volumes, we have the following heat estimates.

\begin{lemma}[{see \cite[Theorem~3.1]{DH1}}]\label{l:heat}
Let $m\geq0$ be an integer and let $\beta,\gamma$ be multi-indices.
For suitable constants $C,c>0$, depending on these orders,
\begin{equation}\label{e:heat}
 |\partial_t^m\partial_x^\beta\partial_y^\gamma h_t(x,y)|
 \leq
 \frac{Ct^{-m-(|\beta|+|\gamma|)/2}}{\mu_\kappa(B(x,\sqrt t))}
 \left(1+\frac{\|x-y\|}{\sqrt t}\right)^{-2}
 e^{-c d(x,y)^2/t}.
\end{equation}

\end{lemma}

Time derivatives of the kernel may be replaced by Dunkl Laplacians, since $\partial_t h_t=\Delta_{\kappa,x}h_t$.

\subsection{Normalization and comparison of averages}\label{s:avg}

After one fixed spatial average is subtracted, the seminorm \eqref{e:bmo} becomes a Hilbert-valued BMO seminorm,
and the resulting function has the growth needed for the heat-kernel integrals.

\begin{lemma}\label{l:norm}
Let $f$ be measurable, with $f(\cdot,t)\in L^1_{\mathrm{loc}}(\mu_\kappa)$
for almost every $t>0$, and suppose that $\Xnorm{f}<\infty$.
Define
$$
 \widetilde f(x,t)=f(x,t)-(f(\cdot,t))_{B(0,1)}.
$$
Then $\widetilde f\in L^2_{\mathrm{loc}}(\mu_\kappa;L^2(\R_+))$ and
$$
 \|\widetilde f\|_{\mathrm{BMO}_2(\mu_\kappa;L^2(\R_+))}
 =\Xnorm{\widetilde f}=\Xnorm{f}.
$$
Moreover, for every $R\geq1$,
\begin{equation}\label{e:ng}
 \left(\frac{1}{\mu_\kappa(B(0,R))}
 \int_{B(0,R)}\|\widetilde f(y,\cdot)\|_{L^2(\R_+)}^2
 \,d\mu_\kappa(y)\right)^{1/2}
 \leq C(1+\log R)\Xnorm{f}.
\end{equation}
\end{lemma}

\begin{proof}
Let $B_j=B(0,2^j)$ for $j\geq0$.
By Cauchy--Schwarz in the spatial variable, we have
$$
 \|(f(\cdot,t))_{B_{j+1}}-(f(\cdot,t))_{B_j}\|_{L^2(\R_+,dt)}
 \leq\left(\frac{\mu_\kappa(B_{j+1})}{\mu_\kappa(B_j)}\right)^{1/2}\Xnorm{f}
 =2^{\mathbf N/2}\Xnorm{f}.
$$
Telescoping from $B_0=B(0,1)$, we obtain
\begin{align*}
 \|\widetilde f\|_{L^2(B_j\times\R_+,d\mu_\kappa\,dt)}
 \leq&\|f(x,t)-(f(\cdot,t))_{B_j}\|_{L^2(B_j\times\R_+,d\mu_\kappa\,dt)}
  +\mu_\kappa(B_j)^{1/2}\|(f(\cdot,t))_{B_j}-(f(\cdot,t))_{B_0}\|_{L^2(\R_+,dt)}\\
 \leq&\mu_\kappa(B_j)^{1/2}
 \bigg(\Xnorm{f}+\sum_{i=0}^{j-1}
 \|(f(\cdot,t))_{B_{i+1}}-(f(\cdot,t))_{B_i}\|_{L^2(\R_+,dt)}\bigg)\\
 \leq&C(1+j)\mu_\kappa(B_j)^{1/2}\Xnorm{f}.
\end{align*}
Choose $j$ with $R\leq2^j<2R$. Doubling then implies \eqref{e:ng}, so
$\widetilde f\in L^2_{\mathrm{loc}}(\mu_\kappa;L^2(\R_+))$.

Since
$$
 \widetilde f(x,t)-(\widetilde f(\cdot,t))_B=f(x,t)-(f(\cdot,t))_B,
$$
we have $\Xnorm{\widetilde f}=\Xnorm{f}$. By Fubini's theorem, we may compute
the $L^2(\R_+)$-valued mean on $B$ by averaging $\widetilde f(\cdot,t)$ in $x$ for almost every $t$. Hence
$$
 \|\widetilde f\|_{\mathrm{BMO}_2(\mu_\kappa;L^2(\R_+))}^2
 =\sup_B\frac{1}{\mu_\kappa(B)}
   \int_0^\infty\int_B
   |\widetilde f(x,t)-(\widetilde f(\cdot,t))_B|^2
   \,d\mu_\kappa(x)\,dt
 =\Xnorm{\widetilde f}^2.
$$
The proof is complete.

\end{proof}

Write $B_*=B(0,1)$ and $a(t)=(f(\cdot,t))_{B_*}$.
Under the assumptions of Lemma~\ref{l:norm}, $\widetilde f$ is integrable
on $B\times(0,T)$ for every ball $B\subset\R^N$ and every $T>0$. Since $f=\widetilde f+a$,
$$
 f\in L^1_{\mathrm{loc}}(\R^N\times[0,\infty),d\mu_\kappa\,dt)
 \quad \text{if and only if}\quad a\in L^1(0,T)\text{ for every }T>0.
$$

For $y\in B(\sigma x,\sqrt t)$, $d(x,y)\leq\|y-\sigma x\|\leq\sqrt t$.
Thus the Gaussian factor $e^{-c d(x,y)^2/t}$ in \eqref{e:heat} is at least $e^{-c}$
even when $\sigma x$ is far from $x$. The next lemma uses the decay in $\|x-y\|$
to control the logarithmic difference of averages on distant balls. For related estimates,
see \cite[Lemma~4.1]{DH3} and \cite[Proposition~4.5]{DH4}.

\begin{lemma}\label{l:mean}
For every $\eta>0$ and $c_0>0$, there is a constant
$C=C(N,R,\kappa,\eta,c_0)$ such that
\[
 \int_{\R^N}
 \frac{e^{-c_0d(x,y)^2/s^2}}{\mu_\kappa(B(x,s))}
 \left(1+\frac{\|x-y\|}{\max(r,s)}\right)^{-\eta}
 \|F(y)-F_{4B}\|_{\mathcal H}\,d\mu_\kappa(y)
  \leq C\bigl(1+|\log(s/r)|\bigr)\Bnorm{F}
\]
for every $F\in\mathrm{BMO}_2(\mu_\kappa;\mathcal H)$, every Euclidean
ball $B=B(c,r)$, every $x\in B$ and every $s>0$.
Here $4B=B(c,4r)$ and $\mathcal H$ may be any separable Hilbert space.
\end{lemma}

The Euclidean decay factor is necessary even when $s=r$, as Proposition~\ref{p:ex} shows.
\begin{proof}
We first compare the averages of $F$ on two balls. If $B\subset B'$,
then by the Cauchy--Schwarz inequality we have
\begin{align*}
 \|F_B-F_{B'}\|_{\mathcal H}
 &=\left\|\frac{1}{\mu_\kappa(B)}
   \int_B(F(y)-F_{B'})\,d\mu_\kappa(y)\right\|_{\mathcal H}\\
 &\leq\left(\frac{1}{\mu_\kappa(B)}
   \int_B\|F(y)-F_{B'}\|_{\mathcal H}^2\,d\mu_\kappa(y)\right)^{1/2}
 \leq\left(\frac{\mu_\kappa(B')}{\mu_\kappa(B)}\right)^{1/2}
   \Bnorm{F}.
\end{align*}
In particular,  by the doubling property of $\mu_{\kappa}$,  we have
$\|F_B-F_{2B}\|_{\mathcal H}\leq C\Bnorm{F}$.

Now let $B=B(c,r)$ and $B'=B(c',r')$ be any two balls, and set
$\rho=r+r'+\|c-c'\|$. Choose $m$ and $n$ so that
$$
 \rho\leq 2^m r<2\rho,\qquad \rho\leq 2^n r'<2\rho.
$$
Both $2^mB$ and $2^nB'$ lie in $B(c,4\rho)$, and by doubling
their measures are comparable to $\mu_\kappa(B(c,4\rho))$.
We compare their averages with $F_{B(c,4\rho)}$ and sum the
differences of averages over successive dilates of $B$ and $B'$.
It follows that
\begin{equation}\label{e:av}
 \|F_B-F_{B'}\|_{\mathcal H}
 \leq C(1+m+n)\Bnorm{F}
 \leq C   \,\bigg(1+\log\frac{\rho^2}{rr'}\,\bigg)\Bnorm{F}.
\end{equation}

Fix $x\in B$ and $s>0$. We apply \eqref{e:av} to
$B(\sigma x,2^js)$ and $4B$, where $\sigma\in G$ and $j\geq0$.
Since $\|\sigma x-c\|\leq\|x-\sigma x\|+r$, we obtain
$$
 \frac{(2^js+4r+\|\sigma x-c\|)^2}{2^js\,4r}
 \leq C\,2^j\frac{\max(r,s)}{\min(r,s)}
 \left(1+\frac{\|x-\sigma x\|}{\max(r,s)}\right)^2.
$$
Furthermore,
\begin{align}\label{e:mean}
 &\frac{1}{\mu_\kappa(B(\sigma x,2^js))}
 \int_{B(\sigma x,2^js)}\|F(y)-F_{4B}\|_{\mathcal H}\,d\mu_\kappa(y)
 \nonumber\\
 \leq\,& \Bnorm{F}+\|F_{B(\sigma x,2^js)}-F_{4B}\|_{\mathcal H}
 \nonumber\\
 \leq \,&  C\left[    1+j+\big|\log(s/r)\big|
 +\log\left(1+\frac{\|x-\sigma x\|}{\max(r,s)}\right)\right]   \Bnorm{F}.
\end{align}

For $y\in B(\sigma x,2^js)$, the triangle inequality gives
$$
 \log\left(1+\frac{\|x-\sigma x\|}{\max(r,s)}\right)
 \leq \log(1+2^j)
     +\log\left(1+\frac{\|x-y\|}{\max(r,s)}\right).
$$
Multiplying by $(1+\|x-y\|/\max(r,s))^{-\eta}$ bounds the right-hand
side by $\log(1+2^j)+1/(e\eta)$, since
$\sup_{u\geq0}(1+u)^{-\eta}\log(1+u)=1/(e\eta)$.
Thus the Euclidean weight controls the last logarithm in
\eqref{e:mean} throughout $B(\sigma x,2^js)$.

Finally, the Gaussian factor satisfies
$$
 e^{-c_0d(x,y)^2/s^2}
 \leq C_{c_0}\sum_{\sigma\in G}\sum_{j=0}^\infty
 e^{-c_0 4^{j-1}}\mathbf1_{B(\sigma x,2^js)}(y).
$$
Indeed, choose $\sigma$ with $\|\sigma x-y\|=d(x,y)$ and take the
smallest $j\geq0$ for which $\|\sigma x-y\|<2^js$.
By $G$-invariance and \eqref{e:vol},
$$
 \mu_\kappa(B(\sigma x,2^js))
 =\mu_\kappa(B(x,2^js))
 \leq C\,2^{j\mathbf N}\mu_\kappa(B(x,s)).
$$
Using \eqref{e:mean} on each ball, we conclude that
\begin{align*}
 &\int_{\R^N}\frac{e^{-c_0d(x,y)^2/s^2}}{\mu_\kappa(B(x,s))}
 \left(1+\frac{\|x-y\|}{\max(r,s)}\right)^{-\eta}
 \|F(y)-F_{4B}\|_{\mathcal H}\,d\mu_\kappa(y)\\
 \leq\,&  C_\eta\Bnorm{F}
 \sum_{\sigma\in G}\sum_{j=0}^\infty
 2^{j\mathbf N}e^{-c_0 4^{j-1}}\left(1+j+|\log(s/r)|\right)\\
 \leq\,& C\left(1+|\log(s/r)|\right)\Bnorm{F}.
\end{align*}
The proof is complete.
\end{proof}

For $t\in[a,b]\subset(0,\infty)$, we use the following lemma to prove absolute convergence
of the spatial integrals involving $h_t$ and its derivatives.

\begin{lemma}\label{l:conv}
Let $\mathcal H_1,\mathcal H_2$ be separable Hilbert spaces and fix
$0<a<b<\infty$. Suppose $K(x,y):\mathcal H_1\to\mathcal H_2$ is
a linear operator such that $(x,y)\mapsto K(x,y)h$ is strongly
measurable for each $h\in\mathcal H_1$ and
\begin{equation}\label{e:fg}
 \|K(x,y)\|\leq C_{a,b}e^{-c d(x,y)^2/b}.
\end{equation}
If $F\in L^2_{\mathrm{loc}}(\mu_\kappa;\mathcal H_1)$ satisfies
$$
 \|F\|_{L^2(B(0,R),\mu_\kappa;\mathcal H_1)}\leq C_F R^m
 \qquad (R\geq1)
$$
for some $C_F,m\geq0$, then
$\int_{\R^N}K(x,y)F(y)\,d\mu_\kappa(y)$ converges absolutely
and defines a locally bounded function of $x$. The integrals over $B(0,R)$
converge locally uniformly in $x$ as $R\to\infty$.

If this kernel represents a bounded operator
$\mathcal T:L^2(\mu_\kappa;\mathcal H_1)\to L^2(\mu_\kappa;\mathcal H_2)$
on bounded compactly supported functions, the same formula holds
almost everywhere for every $F\in L^2(\mu_\kappa;\mathcal H_1)$.
\end{lemma}

\begin{proof}
By \eqref{e:fg} and the Cauchy--Schwarz inequality, for every $R>0$,
$$
 \int_{B(0,R)}\|K(x,y)F(y)\|_{\mathcal H_2}\,d\mu_\kappa(y)
 \leq C_{a,b}\,\mu_\kappa(B(0,R))^{1/2}
       \|F\|_{L^2(B(0,R),\mu_\kappa;\mathcal H_1)}.
$$
For the tail, fix $R_0>0$ and take $R>2R_0+2$. Since $G$ preserves the Euclidean norm,
$$
 d(x,y)\geq\|y\|-\|x\|\geq\|y\|/2
 \qquad (\|x\|\leq R_0,\ \|y\|>R).
$$
On the annulus $2^\ell R<\|y\|\leq 2^{\ell+1}R$, the Gaussian in
\eqref{e:fg} is therefore bounded by $e^{-c4^\ell R^2/(4b)}$.
Applying Cauchy--Schwarz on each annulus, we obtain
\begin{align*}
 \sup_{\|x\|\leq R_0}\int_{\|y\|>R}\|K(x,y)F(y)\|_{\mathcal H_2}\,d\mu_\kappa(y)
 &\leq C_{a,b}\sum_{\ell=0}^\infty e^{-c4^\ell R^2/(4b)}
   \mu_\kappa(B(0,2^{\ell+1}R))^{1/2}
   \|F\|_{L^2(B(0,2^{\ell+1}R),\mu_\kappa;\mathcal H_1)}\\
 &\leq C_{a,b,m}C_F\sum_{\ell=0}^\infty
   (2^\ell R)^{m+\mathbf N/2}e^{-c4^\ell R^2/(4b)}
   \to 0 \qquad (R\to\infty),
\end{align*}
where the last inequality follows from the growth assumption on $F$ and
$\mu_\kappa(B(0,r))=r^{\mathbf N}\mu_\kappa(B(0,1))$.
The two bounds prove absolute convergence and local boundedness. The tails tend to zero
uniformly for $\|x\|\leq R_0$.

Now let $F\in L^2(\mu_\kappa;\mathcal H_1)$ and choose bounded,
compactly supported functions $F_n$ converging to $F$ in this space.
For $F_n-F$, the growth assumption holds with $m=0$ and
$C_F=\|F_n-F\|_{L^2(\mu_\kappa;\mathcal H_1)}$. Taking $R=2R_0+3$
in the two estimates above, we obtain
$$
 \sup_{\|x\|\leq R_0}
 \left\|\int_{\R^N}K(x,y)(F_n(y)-F(y))\,d\mu_\kappa(y)
 \right\|_{\mathcal H_2}
 \leq C_{a,b,R_0}\|F_n-F\|_{L^2(\mu_\kappa;\mathcal H_1)}
 \to 0.
$$
By $L^2$ boundedness, $\mathcal TF_n\to\mathcal TF$ in $L^2(\mu_\kappa;\mathcal H_2)$
and almost everywhere along a subsequence. The kernel formula for $F_{n_k}$ then implies
$$
 \mathcal TF(x)
 =\lim_{k\to\infty}\mathcal TF_{n_k}(x)
 =\lim_{k\to\infty}\int_{\R^N}K(x,y)F_{n_k}(y)\,d\mu_\kappa(y)
 =\int_{\R^N}K(x,y)F(y)\,d\mu_\kappa(y).
$$
The proof is complete.
\end{proof}

For $F\in\mathrm{BMO}_2(\mu_\kappa;\mathcal H)$, the growth assumption of Lemma~\ref{l:conv}
follows from the average comparison \eqref{e:av}:
$$
 \|F\|_{L^2(B(0,R),\mu_\kappa;\mathcal H)}
 \leq C R^{\mathbf N/2}
 \left(\|F_{B(0,1)}\|_{\mathcal H}+(1+\log R)\Bnorm{F}\right),\qquad R\geq1.
$$
Since $\mu_\kappa(B(x,\sqrt t))\geq ct^{\mathbf N/2}$, the bounds in \eqref{e:heat}
imply \eqref{e:fg} on each fixed interval $0<a\leq t\leq b<\infty$.

\subsection{A Hilbert-valued endpoint estimate}\label{s:cz}

The size and first-variable continuity conditions below are the operator-norm versions
of the corresponding conditions in \cite[Definition~1.2]{THp}, with exponent $1/2$; see also \cite{THHLL}.

\begin{proposition}\label{p:cz}
Let $\mathcal H_1,\mathcal H_2$ be separable Hilbert spaces and let
$\mathcal T:L^2(\mu_\kappa;\mathcal H_1)\to L^2(\mu_\kappa;\mathcal H_2)$ be bounded.
Suppose its kernel $K(x,y):\mathcal H_1\to\mathcal H_2$ is strongly measurable, in the sense that $(x,y)\mapsto K(x,y)h$ is strongly measurable for each $h\in\mathcal H_1$,
and satisfies, for $d(x,y)>0$,
\begin{equation}\label{e:ks}
 \|K(x,y)\|
 \leq\frac{C_K}{\mu_\kappa(B(x,d(x,y)))}
 \left(\frac{d(x,y)}{\|x-y\|}\right)^{1/2},
\end{equation}
and
\begin{equation}\label{e:kr}
 \|K(x,y)-K(x',y)\|
 \leq\frac{C_K}{\mu_\kappa(B(x,d(x,y)))}
 \left(\frac{\|x-x'\|}{\|x-y\|}\right)^{1/2},
 \qquad \|x-x'\|\leq d(x,y)/2.
\end{equation}
For $F\in\mathrm{BMO}_2(\mu_\kappa;\mathcal H_1)$, define $\mathcal TF$ by the Bochner integral
$$
 \mathcal TF(x)=\int_{\R^N}K(x,y)F(y)\,d\mu_\kappa(y).
$$
Assume that the integral converges absolutely almost everywhere, agrees with the $L^2$ operator when $F\in L^2$,
and vanishes when $F$ is constant. Then
$$
 \|\mathcal TF\|_{\mathrm{BMO}_2(\mu_\kappa;\mathcal H_2)}
 \leq C\bigl(\|\mathcal T\|_{2\to2}+C_K\bigr)
 \|F\|_{\mathrm{BMO}_2(\mu_\kappa;\mathcal H_1)},
$$
where $C$ is independent of both Hilbert spaces.
\end{proposition}

\begin{proof}
Fix a ball $B=B(c,r)$. We begin by showing that for $x\in B$,
\begin{equation}\label{e:far}
 \int_{(8B)^c}\|K(x,y)-K(c,y)\|\,\|F(y)-F_{8B}\|\,d\mu_\kappa(y)
 \leq CC_K\|F\|_{\mathrm{BMO}_2(\mu_\kappa;\mathcal H_1)}.
\end{equation}
Since $\mathcal T$ vanishes on constants, we may assume by homogeneity that
the BMO seminorm of $F$ is one.

Apply Lemma~\ref{l:mean} to $2B$ with $\eta=1/2$ and $c_0=1$.
Since $4(2B)=8B$ and $e^{-d(z,y)^2/s^2}\geq e^{-1}$ when
$d(z,y)\leq s$, we obtain, for $z\in B$ and $s>0$,
\begin{equation}\label{e:ann}
 \frac{1}{\mu_\kappa(B(z,s))}
 \int_{d(z,y)\leq s}
 \left(1+\frac{\|z-y\|}{\max(2r,s)}\right)^{-1/2}
 \|F(y)-F_{8B}\|\,d\mu_\kappa(y)
 \leq C\left(1+\left|\log\frac{s}{2r}\right|\right).
\end{equation}

First consider $y\notin8B$ with $d(c,y)<2r$. For $z\in B$,
we have $\|z-y\|\geq7r$ and $d(z,y)<3r$. We integrate over
the annuli $s/2<d(z,y)\leq s$, where $s=3r2^{-j}$ and $j\geq0$.
The set $\{y:d(z,y)=0\}$ is finite and hence $\mu_\kappa$-null. On each annulus,
\eqref{e:ks} and doubling imply
\[
 \|K(z,y)\|
  \leq \frac{CC_K}{\mu_\kappa(B(z,s))}
        \left(\frac{s}{\|z-y\|}\right)^{1/2}
 \leq \frac{CC_K}{\mu_\kappa(B(z,s))}
        \left(\frac{s}{r}\right)^{1/2}
        \left(1+\frac{\|z-y\|}{\max(2r,s)}\right)^{-1/2},
\]
where the last inequality uses $\|z-y\|\geq7r$.
Integrating on each annulus and applying \eqref{e:ann}, we obtain
$$
 \int_{\substack{y\notin8B\\d(c,y)<2r}}
 \|K(z,y)\|\,\|F(y)-F_{8B}\|\,d\mu_\kappa(y)
 \leq CC_K\sum_{j=0}^\infty 2^{-j/2}(1+j)
 \leq CC_K.
$$
Taking $z=x$ and $z=c$ bounds the integral of the kernel difference
over this region.

For $d(c,y)\geq2r$, we have $\|x-c\|\leq d(c,y)/2$.
Consider the annulus $s/2\leq d(c,y)<s$, where $s=2^{j+1}r$
and $j\geq1$. Since $\|c-y\|\geq s/2$, applying \eqref{e:kr}
with $c$ as the first point and using doubling, we obtain
\[
 \|K(c,y)-K(x,y)\|
  \leq \frac{CC_K}{\mu_\kappa(B(c,s))}
        \left(\frac{r}{\|c-y\|}\right)^{1/2}
  \leq \frac{CC_K}{\mu_\kappa(B(c,s))}
        \left(\frac{r}{s}\right)^{1/2}
        \left(1+\frac{\|c-y\|}{s}\right)^{-1/2}.
\]
Here $\max(2r,s)=s$, so \eqref{e:ann} with $z=c$ bounds the
contribution of this annulus by $CC_K2^{-j/2}(1+j)$.
Summing over $j\geq1$ completes the proof of \eqref{e:far}.

By hypothesis, the kernel integral for $F-F_{8B}$ converges
absolutely for almost every $x$. Taking one such $x\in B$ in
\eqref{e:far}, we see that the integral defining
$$
 c_B=\int_{(8B)^c}K(c,y)(F(y)-F_{8B})\,d\mu_\kappa(y)
$$
also converges absolutely.

For the integral over $8B$, first restrict $F-F_{8B}$ to the points
of $8B$ where $\|F-F_{8B}\|\leq n$. This truncated function is bounded and
compactly supported, so its kernel integral agrees with the $L^2$
operator. As $n\to\infty$, dominated convergence in $y$ recovers
the integral over $8B$ for almost every $x$. By Fatou's lemma in $x$
and the $L^2$ bound,
\begin{align*}
 \left(\int_B\left\|\int_{8B}K(x,y)(F(y)-F_{8B})\,d\mu_\kappa(y)
    \right\|^2\,d\mu_\kappa(x)\right)^{1/2}
 &\leq\|\mathcal T\|_{2\to2}
   \left(\int_{8B}\|F-F_{8B}\|^2\,d\mu_\kappa\right)^{1/2}\\
 &\leq\|\mathcal T\|_{2\to2}\mu_\kappa(8B)^{1/2}.
\end{align*}

Since $\mathcal T$ vanishes on constants, \eqref{e:far} bounds
the difference between $\mathcal TF(x)-c_B$ and the integral
over $8B$ by $CC_K$ for almost every $x\in B$. Hence
\begin{align*}
 \left(\frac{1}{\mu_\kappa(B)}
 \int_B\|\mathcal TF-(\mathcal TF)_B\|^2\,d\mu_\kappa\right)^{1/2}
 &\leq\left(\frac{1}{\mu_\kappa(B)}
 \int_B\|\mathcal TF-c_B\|^2\,d\mu_\kappa\right)^{1/2}\\
 &\leq\left(\frac{\mu_\kappa(8B)}{\mu_\kappa(B)}\right)^{1/2}
       \|\mathcal T\|_{2\to2}+CC_K\\
 &\leq C\bigl(\|\mathcal T\|_{2\to2}+C_K\bigr).
\end{align*}
Taking the supremum over $B$ completes the proof.
\end{proof}

The following estimate will be used to integrate the heat-kernel bounds in time.

\begin{lemma}\label{l:int}
Fix $c_0>0$ and let $\ell\in\{0,1\}$. For every $x\in\R^N$
and $0<\rho\leq E$, we have
\begin{equation}\label{e:int}
 \left(\int_0^\infty
 \frac{t^{-p\ell/2}e^{-c_0\rho^2/t}}
 {\mu_\kappa(B(x,\sqrt t))^p(1+E/\sqrt t)^p}
 \,\frac{dt}{t}\right)^{1/p}
 \leq \frac{C\rho^{-\ell}}{\mu_\kappa(B(x,\rho))}
 \left(\frac{\rho}{E}\right)^{1/2},
 \qquad p=1,2,
\end{equation}
where $C$ is independent of $x,\rho,E,p$ and $\ell$.
\end{lemma}

\begin{proof}
With $\sqrt t=\rho v$, \eqref{e:vol} implies
$$
 \frac{\mu_\kappa(B(x,\rho))}{\mu_\kappa(B(x,\rho v))}
 \leq\begin{cases}
 Cv^{-\mathbf N},&0<v\leq1,\\[4pt]
 Cv^{-N},&v\geq1,
 \end{cases}
$$
and $\left(1+\frac E{\rho v}\right)^{-1}
 \leq\left(\frac\rho E\right)^{1/2}v^{1/2}.$
Therefore the left-hand side of \eqref{e:int} is at most
\[
 \frac{C\rho^{-\ell}}{\mu_\kappa(B(x,\rho))}
 \left(\frac\rho E\right)^{1/2}
 \left(
 \int_0^1v^{p(1/2-\ell-\mathbf N)}e^{-c_0/v^2}\frac{dv}{v}
 +\int_1^\infty v^{p(1/2-\ell-N)}\frac{dv}{v}
 \right)^{1/p}
 \leq \frac{C\rho^{-\ell}}{\mu_\kappa(B(x,\rho))}
 \left(\frac\rho E\right)^{1/2},
\]
due to  $N+\ell>1/2$.
The proof is complete.
\end{proof}

To apply Proposition~\ref{p:cz} in Sections~\ref{sec:3} and \ref{sec:4},
we use Lemma~\ref{l:int} to show that each operator-valued kernel $K$ satisfies
$$
 \|\partial_x^\beta K(x,y)\|
 \leq \frac{Cd(x,y)^{-|\beta|}}{\mu_\kappa(B(x,d(x,y)))}
 \left(\frac{d(x,y)}{\|x-y\|}\right)^{1/2},
 \qquad d(x,y)>0,\quad |\beta|\leq1.
$$
When $\beta=0$, this is the size estimate \eqref{e:ks}.
To obtain \eqref{e:kr}, fix $x'$ with $\|x-x'\|\leq d(x,y)/2$.
For $z$ on the segment from $x$ to $x'$, the triangle inequality implies
$d(z,y)\simeq d(x,y)$ and $\|z-y\|\simeq\|x-y\|$. By doubling,
$\mu_\kappa(B(z,d(z,y)))\simeq\mu_\kappa(B(x,d(x,y)))$.
Integrating the bound for $|\beta|=1$ along this segment, we obtain
\begin{equation}\label{e:kd}
\begin{split}
 \|K(x,y)-K(x',y)\|
 &\leq \frac{C}{\mu_\kappa(B(x,d(x,y)))}
 \frac{\|x-x'\|}{d(x,y)}
 \left(\frac{d(x,y)}{\|x-y\|}\right)^{1/2}\\
 &\leq \frac{C}{\mu_\kappa(B(x,d(x,y)))}
 \left(\frac{\|x-x'\|}{\|x-y\|}\right)^{1/2}.
\end{split}
\end{equation}

\subsection{The spatial and time norms}\label{s:ex}

The following rank-one example shows that the Gaussian factor $e^{-c_0d(x,y)^2/s^2}$ alone does not suffice in Lemma~\ref{l:mean}.

\begin{proposition}\label{p:ex}
Let $R=\{-1,1\}$ and $\kappa\geq0$, so that $d\mu_\kappa(x)=|x|^{2\kappa}\,dx$ on $\R$.
The function $b(x)=\log|x-1|$ belongs to $\mathrm{BMO}_2(\mu_\kappa)$.
For $r>0$, let $I_r=(1-r,1+r)$ and denote its reflection by $-I_r=(-1-r,-1+r)$. Then
\begin{equation}\label{e:ex}
 b_{I_r}=\log r-1+o(1),\qquad b_{-I_r}=\log2+o(1),\qquad r\to0^+.
\end{equation}
For every $c_0>0$ and $0<r<1/8$,
\begin{equation}\label{e:fail}
 \frac1{\mu_\kappa(I_r)}\int_\R
 e^{-c_0d(1,y)^2/r^2}|b(y)-b_{4I_r}|\,d\mu_\kappa(y)
 \geq c\log(1/r)-C,
\end{equation}
where $4I_r=I_{4r}$ and $c,C>0$ depend only on $\kappa$ and $c_0$.
\end{proposition}

\begin{proof}
Let $I=(z-r,z+r)$. If $|z-1|\geq2r$, then
$\sup_I b-\inf_I b\leq\log3$.
Otherwise, subtracting $\log r$ and using
$\mu_\kappa(I)\simeq r(|z|+r)^{2\kappa}$, we obtain
\[
 \frac1{\mu_\kappa(I)}\int_I|b(x)-b_I|^2\,d\mu_\kappa(x)
 \leq \frac Cr\int_I\left|\log\frac{|x-1|}{r}\right|^2\,dx
 \leq C\int_{-3}^3|\log|v||^2\,dv\leq C.
\]
Thus $b\in\mathrm{BMO}_2(\mu_\kappa)$.

With $x=1+rv$,
$$
 b_{I_r}-\log r
 =\frac{\int_{-1}^1\log|v|\,(1+rv)^{2\kappa}\,dv}
 {\int_{-1}^1(1+rv)^{2\kappa}\,dv}
 \to\frac12\int_{-1}^1\log|v|\,dv=-1.
$$
On $-I_r$, $b=\log2+o(1)$ uniformly as $r\to0+$.
The two disjoint intervals have equal measure by reflection invariance, so
$b_{I_r\cup(-I_r)}=(b_{I_r}+b_{-I_r})/2$.
Jensen's inequality on each interval then implies
$$
 \frac1{\mu_\kappa(I_r\cup(-I_r))}
 \int_{I_r\cup(-I_r)}|b-b_{I_r\cup(-I_r)}|\,d\mu_\kappa
 \geq\frac12|b_{I_r}-b_{-I_r}|
 =\frac12\log(1/r)+O(1).
$$
Thus $b$ does not belong to orbit BMO (see \cite[Section~7.2]{JL} for the strict inclusion in Euclidean BMO).
For $y\in-I_r$, we have $d(1,y)=|1+y|\leq r$.
Restricting the integral on the left of \eqref{e:fail} to $-I_r$ and applying Jensen's inequality, we obtain
\[
 \frac1{\mu_\kappa(I_r)}\int_\R
 e^{-c_0d(1,y)^2/r^2}|b(y)-b_{4I_r}|\,d\mu_\kappa(y)
 \geq e^{-c_0}|b_{-I_r}-b_{4I_r}|
 = e^{-c_0}\log(1/r)+O(1),
\]
where we used \eqref{e:ex} with radii $r$ and $4r$.
The proof is complete.
\end{proof}

Since $|1-y|\simeq1$ on $-I_r$, for every fixed $\eta>0$,
$$
 \frac1{\mu_\kappa(I_r)}\int_{-I_r}
 \frac{e^{-c_0d(1,y)^2/r^2}}{(1+|1-y|/r)^\eta}
 |b(y)-b_{4I_r}|\,d\mu_\kappa(y)
 \simeq_{\kappa,c_0,\eta} r^\eta\log(1/r)
 \qquad(r\to0+).
$$
The logarithm in \eqref{e:mean} cannot be replaced by a constant independent of $\|x-\sigma x\|$.

For instance, Theorem~\ref{t:main} applies to
$$
 f(x,t)=\chi(t)\log|x-1|,\qquad 0\ne\chi\in C_c^\infty(\R_+).
$$
Then
$$
 \Xnorm{f}=\|\chi\|_{L^2(\R_+)}\|b\|_{\mathrm{BMO}_2(\mu_\kappa)}<\infty,
$$
whereas the seminorm obtained by using orbit intervals in \eqref{e:bmo} is infinite.

The maximal regularity estimate fails when orbit intervals are used on both sides.
For $d\mu_\kappa(x)=|x|^{2\kappa}\,dx$ on $\R$, $\kappa\geq0$, write
\begin{equation}\label{e:orb}
 \|g\|_{\mathrm{orb}}^2
 =\sup_I\frac1{\mu_\kappa(I\cup(-I))}
 \int_0^\infty\int_{I\cup(-I)}
 |g(x,t)-(g(\cdot,t))_{I\cup(-I)}|^2\,d\mu_\kappa(x)\,dt,
\end{equation}
where $I$ ranges over bounded intervals of positive length, $-I=\{-x:x\in I\}$,
and the subscript denotes the $\mu_\kappa$-average over $I\cup(-I)$.

\begin{proposition}\label{p:orb}
For every $\kappa\geq0$ and every integer $m\geq1$, there is a bounded external force $f_m$,
odd in $x$ and compactly supported in $\R\times(0,1)$, whose Duhamel solution $u_m$ in \eqref{e:duh} satisfies
$$
 \|f_m\|_{\mathrm{orb}}\leq1,
 \qquad
 \|\Delta_\kappa u_m\|_{\mathrm{orb}}\geq c_\kappa\sqrt m,
$$
where $c_\kappa>0$ is independent of $m$.
\end{proposition}

\begin{proof}
For an odd function $g(\cdot,t)$, the average over $I\cup(-I)$ vanishes.
Weighted Lebesgue differentiation, first over a finite time interval, implies
\begin{equation}\label{e:oid}
 \|g\|_{\mathrm{orb}}^2
 =\operatorname*{ess\,sup}_{x>0}\int_0^\infty|g(x,t)|^2\,dt.
\end{equation}
We use disjoint spatial supports to bound the external force in \eqref{e:oid}
and estimate $\Delta_\kappa u_m(1,\cdot)$ on disjoint time intervals.
Fix $0\not\equiv\varphi\in C_c^\infty((1,2))$ with $0\leq\varphi\leq1$.
For $0<\varepsilon<1/4$, to be fixed below, let
$$
 r_j=\varepsilon4^{-j},\qquad t_j=1-2^{-j},\qquad j\geq1,
$$
and define
$$
 \varphi_j(x)=\varphi\left(\frac{x-1}{r_j}\right)
 -\varphi\left(\frac{-x-1}{r_j}\right),
 \qquad
 f_m(x,t)=\sum_{j=1}^m r_j^{-1}\varphi_j(x)\mathbf1_{(t_j,t_j+r_j^2)}(t).
$$
The functions $\varphi_j$ are odd and have pairwise disjoint supports, and
the time intervals are contained in $(0,1)$. Thus \eqref{e:oid} implies
$$
 \|f_m\|_{\mathrm{orb}}^2
 =\operatorname*{ess\,sup}_{x>0}\sum_{j=1}^m|\varphi_j(x)|^2\leq1.
$$

The heat semigroup preserves oddness, so both $u_m$ and $\Delta_\kappa u_m$ are odd.
 To estimate $(e^{s\Delta_\kappa}\varphi_j)(1)$, take $s=ar_j^2$ with $1/2\leq a\leq1$.
For $y=1+r_jz$, $1\leq z\leq2$, we have $y\simeq1$ and
$(1-y)^2/(4s)=z^2/(4a)\leq2$. If $\kappa>0$,
the two estimates in \cite[Proposition~2.3(e)]{ABDH} at $(1,y)$ and $(1,-y)$ imply
$$
 h_s(1,y)\geq c_\kappa r_j^{-1},
 \qquad h_s(1,-y)\leq C_\kappa r_j.
$$
For $\kappa=0$, the first bound follows from the Gaussian formula, and
$h_s(1,-y)\leq C r_j^{-1}e^{-1/r_j^2}\leq C r_j$.
Choose $\varepsilon$ small enough that
$h_s(1,y)-h_s(1,-y)\geq c r_j^{-1}$ throughout this range.
By oddness of $\varphi_j$,
\begin{equation}\label{e:olb}
 (e^{a r_j^2\Delta_\kappa}\varphi_j)(1)
 = \int_{1+r_j}^{1+2r_j}
 [h_s(1,y)-h_s(1,-y)]\,
 \varphi\left(\frac{y-1}{r_j}\right)y^{2\kappa}\,dy
 \geq c\int_1^2\varphi(z)\,dz=:c_\kappa>0.
\end{equation}
For the terms with $i<j$, \eqref{e:heat}, \eqref{e:vol}, and
$\|\varphi_i\|_{L^1(\mu_\kappa)}\leq C_\kappa r_i$ imply
\begin{equation}\label{e:od}
 \left|\partial_s(e^{s\Delta_\kappa}\varphi_i)(1)\right|
 \leq\frac{C_\kappa\|\varphi_i\|_{L^1(\mu_\kappa)}}{s\mu_\kappa(B(1,\sqrt s))}
 \leq C_\kappa r_i s^{-3/2},\qquad 0<s\leq1.
\end{equation}

For $1\leq j\leq m$, let $J_j=(t_j+r_j^2/2,t_j+r_j^2)$.
Since $t_j+r_j^2<t_{j+1}$, the intervals $J_j$ are pairwise disjoint.
For $t\in J_j$, only the terms with $i\leq j$ contribute to \eqref{e:duh}.
Using $\partial_s e^{s\Delta_\kappa}=\Delta_\kappa e^{s\Delta_\kappa}$ and $\varphi_i(1)=0$, we obtain
\[
 \Delta_\kappa u_m(1,t)
 = r_j^{-1}(e^{(t-t_j)\Delta_\kappa}\varphi_j)(1)
 +\sum_{i<j}r_i^{-1}
 \bigl[(e^{(t-t_i)\Delta_\kappa}-e^{(t-t_i-r_i^2)\Delta_\kappa})\varphi_i\bigr](1).
\]
For $i<j$, we have $t-t_i-r_i^2\geq2^{-i-2}$ and $t-t_i<1$.
Consequently, \eqref{e:od} implies
\begin{align*}
 \left|\sum_{i<j}r_i^{-1}
 \bigl[(e^{(t-t_i)\Delta_\kappa}-e^{(t-t_i-r_i^2)\Delta_\kappa})\varphi_i\bigr](1)\right|
 \leq&C_\kappa\sum_{i<j}\int_{t-t_i-r_i^2}^{t-t_i}s^{-3/2}\,ds\\
 \leq&C_\kappa\sum_{i<j}r_i^2\,2^{3i/2}
 \leq C_\kappa\varepsilon^2\sum_{i\geq1}2^{-5i/2}
 \leq C_\kappa\varepsilon^2.
\end{align*}
Decrease $\varepsilon$ so that $C_\kappa\varepsilon^2\leq c_\kappa/(2r_1)$.
Since $r_j\leq r_1$, \eqref{e:olb} now implies
$$
 |\Delta_\kappa u_m(1,t)|\geq\frac{c_\kappa}{2r_j},\qquad t\in J_j.
$$
For fixed $m$, all $\varphi_i$, $i\leq m$, vanish on $|x-1|<r_m/2$.
The same finite-sum formula is valid there and is continuous in $x$,
uniformly on the closures of $J_1,\ldots,J_m$.
Using oddness, \eqref{e:oid}, and $|J_j|=r_j^2/2$, we conclude that
\[
 \|\Delta_\kappa u_m\|_{\mathrm{orb}}^2
 \geq \sum_{j=1}^m\int_{J_j}|\Delta_\kappa u_m(1,t)|^2\,dt
 \geq\frac{c_\kappa^2}{4}\sum_{j=1}^m\frac{|J_j|}{r_j^2}
 =\frac{c_\kappa^2m}{8}.
\]
The proof is complete.
\end{proof}

\begin{remark}\label{r:time}
For each ball $B$ and each $t$ with $f(\cdot,t)\in L^2(B,\mu_\kappa)$,
$$
 \frac1{\mu_\kappa(B)^2}\int_{B\times B}|f(x,t)-f(y,t)|^2\,d\mu_\kappa(x)\,d\mu_\kappa(y)
 =\frac2{\mu_\kappa(B)}\int_B|f(x,t)-(f(\cdot,t))_B|^2\,d\mu_\kappa(x).
$$
When $\kappa=0$, the square of the Ogawa--Shimizu seminorm in \cite[(1.10)]{OS} with $\theta=2$
is twice the right-hand side of \eqref{e:bmo} with time integral over $(0,r_B^2)$.
For every $\kappa\geq0$, this shorter interval defines a strictly larger space.
Indeed, take $0\ne\varphi\in C_c^\infty(\R^N)$ and
$f(x,t)=\mathbf1_{(1,\infty)}(t)t^{-1/2}\varphi(x)$.
The integral over $(0,r_B^2)$ vanishes when $r_B\leq1$. For $r_B>1$, \eqref{e:vol} implies
$$
 \frac1{\mu_\kappa(B)}\int_0^{r_B^2}\int_B
 |f(x,t)-(f(\cdot,t))_B|^2\,d\mu_\kappa(x)\,dt
 \leq\frac{2\log r_B}{\mu_\kappa(B)}\|\varphi\|_{L^2(\mu_\kappa)}^2
 \leq C\frac{\log r_B}{r_B^N}\|\varphi\|_{L^2(\mu_\kappa)}^2.
$$
This is uniformly bounded. On any ball where $\varphi$ is not constant, the integral over $(0,\infty)$ is infinite,
since its positive spatial variance is multiplied by $\int_1^\infty dt/t$.
\end{remark}

In the Bochner norm the ball may vary with time. We use a different ball on each time interval
and place the cosine oscillations in balls far from the reflecting hyperplanes.

\begin{proposition}\label{p:bm}
Let $0<T\leq\infty$.
There is no constant $C_T$ such that
$$
 \|\Delta_\kappa u\|_{L^2(0,T;\mathrm{BMO}(\mu_\kappa))}
 \leq C_T\|f\|_{L^2(0,T;\mathrm{BMO}(\mu_\kappa))}
$$
for every $f\in C_c^\infty(\R^N\times(0,T))$, where $u$ is its Duhamel solution \eqref{e:duh}.
The estimate also fails with $\partial_tu$ in place of $\Delta_\kappa u$.
\end{proposition}

\begin{proof}
We first treat $T=3$.
We construct uniformly bounded external forces $f_m$ whose solutions satisfy
$\|\partial_tu_m(\cdot,t)\|_{\mathrm{BMO}(\mu_\kappa)}\geq c\sqrt m$ for every $1\leq t\leq2$.
All constants below are independent of $m$.

Fix $r=2\pi\sqrt N$ and choose $\chi\in C_c^\infty(B(0,2r))$ with $0\leq\chi\leq1$
and $\chi=1$ on $B(0,r)$.
Choose a unit vector $\theta$ outside the reflecting hyperplanes, $S>4r$, and $A>0$ such that
$A|\langle\alpha,\theta\rangle|>4r\|\alpha\|$ for every $\alpha\in R_+$.
Write
$$
 x_k=(A+kS)\theta,\qquad B_k=B(x_k,r),\qquad
 Q_k=x_k+(-\pi,\pi)^N,\qquad k\geq0.
$$
The balls $B(x_k,2r)$ are disjoint and have distance at least $2r$ from every reflecting hyperplane.
Moreover, $Q_k\subset B_k$ and, uniformly in $k$,
\begin{equation}\label{e:bw}
 d\mu_\kappa(x)\simeq
 \left(\prod_{\alpha\in R_+}|\langle\alpha,x_k\rangle|^{2\kappa(\alpha)}\right)dx,
 \qquad x\in B_k.
\end{equation}

Fix a nonzero function $\rho\in C_c^\infty((1,2))$ with $0\leq\rho\leq1$.
For an integer $m\geq1$, let $\delta_j=4^{-j}$ for $1\leq j\leq m$, and let
$$
 t_k=1+k\delta_m,\qquad I_k=[t_k,t_k+\delta_m],\qquad 0\leq k<4^m.
$$
For these indices define
$$
 a_{kj}(t)=\int_0^t e^{-(t-s)/\delta_j}
 \rho\left(\frac{t_k-s}{\delta_j}\right)\,ds,\qquad t\geq0.
$$
Then
$$
 a'_{kj}(t)+\delta_j^{-1}a_{kj}(t)
 =\rho\left(\frac{t_k-t}{\delta_j}\right),\qquad
 a_{kj}(0)=0,\qquad 0\leq a_{kj}(t)\leq\delta_j.
$$
For each $k$, the intervals $(t_k-2\delta_j,t_k-\delta_j)$, $1\leq j\leq m$,
are pairwise disjoint and lie in $(1/2,2)$. Hence
$$
 \sup_{k,t\geq0}\sum_{j=1}^m\rho\left(\frac{t_k-t}{\delta_j}\right)\leq1.
$$

Choose $\eta\in C_c^\infty((0,3))$ with $\eta=1$ on $[1/2,5/2]$, and define
$$
 \Phi_k(x,t)=\sum_{j=1}^m a_{kj}(t)\cos\bigl(2^j(x_1-x_{k,1})\bigr),\qquad
 u_m(x,t)=\eta(t)\sum_{k=0}^{4^m-1}\chi(x-x_k)\Phi_k(x,t).
$$
By the defining equations for $a_{kj}$,
$$
 (\partial_t-\Delta)\Phi_k(x,t)
 =\sum_{j=1}^m\rho\left(\frac{t_k-t}{\delta_j}\right)
 \cos\bigl(2^j(x_1-x_{k,1})\bigr).
$$
The spatial cutoffs have disjoint supports. The product rule therefore yields
$$
 \|u_m\|_\infty+\|\nabla u_m\|_\infty
 +\|(\partial_t-\Delta)u_m\|_\infty
 \leq C\left(1+\sum_{j=1}^m(1+2^j)\delta_j\right)\leq C.
$$
The support of $u_m$ stays a fixed distance from the reflecting hyperplanes.
When the reflection difference is nonzero, $x$ or $\sigma_\alpha x$ is in this support.
Since $|\langle\alpha,\sigma_\alpha x\rangle|=|\langle\alpha,x\rangle|$,
the denominators in \eqref{e:lap} are bounded away from zero on all nonzero terms. Hence
$$
 \|(\Delta_\kappa-\Delta)u_m\|_\infty
 \leq C\bigl(\|\nabla u_m\|_\infty+\|u_m\|_\infty\bigr)\leq C.
$$
The spatial support of $(\Delta_\kappa-\Delta)u_m$ is contained in the union of finitely many reflections of that of $u_m$.
Thus $f_m=(\partial_t-\Delta_\kappa)u_m$ belongs to $C_c^\infty(\R^N\times(0,3))$, and
\begin{equation}\label{e:bf}
 \|f_m\|_\infty+\|f_m\|_{L^2(0,3;\mathrm{BMO}(\mu_\kappa))}\leq C.
\end{equation}
For each $m$, the functions $f_m$, $\partial_tu_m$ and $\Delta_\kappa u_m$ are continuous
as $L^\infty$-valued functions of $t$, hence strongly measurable in BMO.
Since $u_m(\cdot,0)=0$, integration in $L^2(\mu_\kappa)$ yields the Duhamel formula
\[
 \int_0^t e^{(t-s)\Delta_\kappa}f_m(\cdot,s)\,ds
 = \int_0^t\frac{d}{ds}\left(e^{(t-s)\Delta_\kappa}u_m(\cdot,s)\right)\,ds
 =u_m(\cdot,t),\qquad t>0.
\]
Here $u_m$ and $f_m$ are extended by zero outside $(0,3)$.

Let $c_\rho=\int_1^2e^{-s}\rho(s)\,ds>0$.
For $t\in I_k$, all the source terms in the definition of $a_{kj}$ vanish. Changing variables in the definition of $a_{kj}$, we have
$$
 a_{kj}(t)=c_\rho\delta_j e^{-(t-t_k)/\delta_j},\qquad 1\leq j\leq m.
$$
On $Q_k$, the spatial cutoff and $\eta$ are both equal to one. Therefore,
$$
 \partial_tu_m(x,t)
 =-c_\rho\sum_{j=1}^m e^{-(t-t_k)/\delta_j}
 \cos\bigl(2^j(x_1-x_{k,1})\bigr),\qquad x\in Q_k,\quad t\in I_k.
$$
The cosines are orthogonal and have mean zero on $Q_k$.
Using \eqref{e:bw} and $0\leq t-t_k\leq\delta_m\leq\delta_j$, we obtain
\begin{align*}
 \|\partial_tu_m(\cdot,t)\|_{\mathrm{BMO}(\mu_\kappa)}^2
 \geq&\inf_{b\in\mathbb C}\frac1{\mu_\kappa(B_k)}
 \int_{B_k}|\partial_tu_m(x,t)-b|^2\,d\mu_\kappa(x)\\
 \geq&c\inf_{b\in\mathbb C}\frac1{|Q_k|}
 \int_{Q_k}|\partial_tu_m(x,t)-b|^2\,dx\\
 =&\frac{cc_\rho^2}{2}\sum_{j=1}^m e^{-2(t-t_k)/\delta_j}
 \geq\frac{cc_\rho^2e^{-2}}2\,m.
\end{align*}
The intervals $I_k$ cover $[1,2]$.
Together with \eqref{e:bf} and $\Delta_\kappa u_m=\partial_tu_m-f_m$, this proves
$$
 \|\partial_tu_m\|_{L^2(0,3;\mathrm{BMO}(\mu_\kappa))}\geq  c\sqrt m \qquad \text{and}\qquad
 \|\Delta_\kappa u_m\|_{L^2(0,3;\mathrm{BMO}(\mu_\kappa))}\geq  c\sqrt m-C.
$$
Thus neither maximal regularity estimate holds for $T=3$ or, by zero extension, for $T=\infty$.
For a general finite interval $(0,T)$, set $a=T/3$.
The rescaled solution and external force
$$
 a\,u_m\left(\frac{x}{\sqrt a},\frac{t}{a}\right),
 \qquad
 f_m\left(\frac{x}{\sqrt a},\frac{t}{a}\right)
$$
satisfy \eqref{e:pde} with zero initial value on $\R^N\times(0,T)$.
Both time norms are multiplied by $a^{1/2}$, since the spatial BMO seminorm is dilation-invariant.
The proof is complete.
\end{proof}

\section{Maximal regularity and the Cauchy problem}\label{sec:3}

\subsection{The operator estimate}\label{s:mr}

For $0<a<b<\infty$, define
\begin{equation}\label{e:cut}
 \mathcal M_\kappa^{a,b}f(x,t)
 =\int_a^b\Delta_\kappa e^{\tau\Delta_\kappa}f(\cdot,t-\tau)(x)\,d\tau.
\end{equation}
Here we set $f(x,s)=0$ for $s<0$, so the integrand vanishes when $\tau>t$.

\begin{proposition}\label{p:mr}
For every $f\in L^2_{\mathrm{loc}}(\mu_\kappa;L^2(\R_+))$ with $\Xnorm{f}<\infty$,
$$
 \Xnorm{\mathcal M_\kappa^{a,b}f}\leq C\Xnorm{f},\qquad 0<a<b<\infty,
$$
where $C$ is independent of $a,b$ and $f$.
\end{proposition}

\begin{proof}
For fixed $a,b$, Lemmas~\ref{l:heat} and \ref{l:conv} ensure
$\mathcal M_\kappa^{a,b}f\in L^2_{\mathrm{loc}}(\mu_\kappa;L^2(\R_+))$.
For $f\in L^2(\mu_\kappa\times\R_+)$, we take the Fourier transform in time
and the spectral representation of $-\Delta_\kappa$ in space. The multiplier of \eqref{e:cut} is
\begin{equation}\label{e:mult}
 -\int_a^b\lambda e^{-\tau(\lambda+i\xi)}\,d\tau
 =-\frac\lambda{\lambda+i\xi}
 \left(e^{-a(\lambda+i\xi)}-e^{-b(\lambda+i\xi)}\right),\qquad\lambda\geq0.
\end{equation}
The multiplier is bounded in modulus by $2$. As $a\to0+$ and $b\to\infty$,
it tends to $-\lambda/(\lambda+i\xi)$ for $\lambda>0$ and is zero at $\lambda=0$. Hence
\begin{equation}\label{e:l2}
 \|\mathcal M_\kappa^{a,b}f\|_{L^2(\mu_\kappa\times\R_+)}
 \leq2\|f\|_{L^2(\mu_\kappa\times\R_+)}.
\end{equation}

In Proposition~\ref{p:cz}, we use $\mathcal H_1=\mathcal H_2=L^2(\R_+)$.
The spatial kernel acts by
$$
 K^{a,b}(x,y)g
 =\int_a^b\partial_\tau h_\tau(x,y)g(\cdot-\tau)\,d\tau,
 \qquad g\in L^2(\R_+),
$$
where the integral is taken strongly in $L^2(\R_+)$ and $g$ is extended by zero for negative times.
Time shifts are contractions in this space. Using \eqref{e:heat} and Lemma~\ref{l:int} with $p=1$,
we obtain, for $|\beta|\leq1$ and $d(x,y)>0$,
\begin{align*}
 \|\partial_x^\beta K^{a,b}(x,y)\|
 \leq&\int_a^b|\partial_x^\beta\partial_\tau h_\tau(x,y)|\,d\tau\\
 \leq&C\int_0^\infty
 \frac{\tau^{-|\beta|/2}e^{-cd(x,y)^2/\tau}}
 {\mu_\kappa(B(x,\sqrt\tau))(1+\|x-y\|/\sqrt\tau)}\frac{d\tau}{\tau}\\
 \leq&\frac{Cd(x,y)^{-|\beta|}}{\mu_\kappa(B(x,d(x,y)))}
 \left(\frac{d(x,y)}{\|x-y\|}\right)^{1/2}.
\end{align*}
These estimates imply \eqref{e:ks} and, by \eqref{e:kd}, \eqref{e:kr}, uniformly in $a,b$.

With $a,b$ fixed, \eqref{e:heat} and the lower volume bound
$\mu_\kappa(B(x,\sqrt\tau))\geq c\tau^{\mathbf N/2}$ also show that
$$
 \|K^{a,b}(x,y)\|
 \leq\int_a^b|\partial_\tau h_\tau(x,y)|\,d\tau
 \leq C_{a,b}e^{-c d(x,y)^2/b}.
$$
For $f$ as in Proposition~\ref{p:mr}, the Bochner integral converges absolutely by this bound and Lemma~\ref{l:conv}.
For bounded compactly supported $L^2(\R_+)$-valued functions, the heat-kernel formula
and Fubini's theorem show that the kernel integral coincides with \eqref{e:cut}.
By Lemma~\ref{l:conv} and \eqref{e:l2}, the two agree for all $f\in L^2(\mu_\kappa;L^2(\R_+))$.
For a spatial constant $g\in L^2(\R_+)$, absolute convergence and the conservation property $\int_{\R^N}h_\tau(x,y)\,d\mu_\kappa(y)=1$ yield
$$
 \int_{\R^N}K^{a,b}(x,y)g\,d\mu_\kappa(y)
 =\int_a^b g(\cdot-\tau)
   \left(\int_{\R^N}\partial_\tau h_\tau(x,y)\,d\mu_\kappa(y)\right)d\tau=0.
$$
Proposition~\ref{p:cz} now implies
$$
 \Xnorm{\mathcal M_\kappa^{a,b}f}
 =\|\mathcal M_\kappa^{a,b}f\|_{\mathrm{BMO}_2(\mu_\kappa;L^2(\R_+))}
 \leq C\|f\|_{\mathrm{BMO}_2(\mu_\kappa;L^2(\R_+))}
 =C\Xnorm{f}.
$$
The proof is complete.
\end{proof}

\subsection{Construction and uniqueness}\label{s:sol}

To pass from oscillation bounds to local $L^2$ estimates, we need the spatial mean of the Laplacian.

\begin{lemma}\label{l:rec}
Let $0<T<\infty$ and $v\in C([0,T];L^2_{\mathrm{loc}}(\mu_\kappa))$.
Suppose $\Delta_\kappa v=g$ in weighted distributions on $\R^N\times(0,T)$,
where $g\in L^1_{\mathrm{loc}}(\R^N\times(0,T),d\mu_\kappa\,dt)$ and
$$
 \sup_B\frac1{\mu_\kappa(B)}
 \int_0^T\int_B|g(x,t)-(g(\cdot,t))_B|^2\,d\mu_\kappa(x)\,dt<\infty.
$$
Then $g\in L^2(B\times(0,T),d\mu_\kappa\,dt)$ for every ball $B$.
\end{lemma}

\begin{proof}
Choose a radial $\phi\in C_c^\infty(B_*)$ with $\int\phi\,d\mu_\kappa=1$.
Then $\Delta_\kappa\phi$ is also supported in $B_*$.
Testing the distributional identity against $\phi(x)\eta(t)$, where $\eta\in C_c^\infty(0,T)$,
we obtain, for almost every $t\in(0,T)$,
$$
 (g(\cdot,t))_{B_*}
 =\int_{B_*}v(x,t)\Delta_\kappa\phi(x)\,d\mu_\kappa(x)
 -\int_{B_*}\bigl(g(x,t)-(g(\cdot,t))_{B_*}\bigr)\phi(x)\,d\mu_\kappa(x).
$$
Consequently,
\begin{align*}
 \|(g(\cdot,t))_{B_*}\|_{L^2(0,T)}
 &\leq\sqrt T\,\|\Delta_\kappa\phi\|_{L^2(\mu_\kappa)}
 \sup_{0\leq t\leq T}\|v(t)\|_{L^2(B_*,\mu_\kappa)}\\
 &\quad+\|\phi\|_{L^2(\mu_\kappa)}
 \left(\int_0^T\int_{B_*}|g(x,t)-(g(\cdot,t))_{B_*}|^2\,d\mu_\kappa(x)\,dt\right)^{1/2}<\infty.
\end{align*}
Apply Lemma~\ref{l:norm} to $g$ extended by zero outside $(0,T)$.
On each $B\times(0,T)$ the difference $g-(g(\cdot,t))_{B_*}$ is square-integrable,
and the mean belongs to $L^2(0,T)$.
The proof is complete.
\end{proof}

We next estimate the spatial tails of the heat integrals for $\widetilde f$.

\begin{lemma}\label{l:tail}
Let $f$ and $\widetilde f$ be as in Lemma~\ref{l:norm}, and fix $R_0,T>0$.
For each of the kernels $q_\tau=h_\tau$ and $q_\tau=\Delta_{\kappa,x}h_\tau$, the integral
$$
 \int_0^t\int_{\|y\|>R}|q_{t-s}(x,y)\widetilde f(y,s)|\,d\mu_\kappa(y)\,ds
$$
belongs to $L^2(0,T)$ for $\|x\|\le R_0$ and $R\ge2\max\{1,R_0,\sqrt T\}$.
Its $L^2(0,T)$ norm tends to zero uniformly for $\|x\|\le R_0$ as $R\to\infty$.
For $q_\tau=h_\tau$, the integral also tends to zero in $L^\infty(0,T;L^2(B(0,R_0),\mu_\kappa))$.
Both limits hold uniformly if $t-s$ is restricted to any measurable subset of $(0,T)$.
\end{lemma}

\begin{proof}
Set $\ell=0$ for $q_\tau=h_\tau$ and $\ell=1$ for
$q_\tau=\Delta_{\kappa,x}h_\tau=\partial_\tau h_\tau$.
Estimate \eqref{e:heat}, with $m=\ell$ and $\beta=\gamma=0$, implies
$$
 |q_\tau(x,y)|\leq C\tau^{-\ell}\mu_\kappa(B(x,\sqrt\tau))^{-1}e^{-c d(x,y)^2/\tau}.
$$
By \eqref{e:vol}, $\mu_\kappa(B(x,\sqrt\tau))\geq c\tau^{\mathbf N/2}$.
For $\|x\|\leq R_0$ and $\|y\|>R$, we have $d(x,y)\geq\|y\|-\|x\|\geq\|y\|/2$.
Let $R_j=2^jR$ and $E_j=B(0,2R_j)\setminus B(0,R_j)$.
Cauchy--Schwarz and \eqref{e:ng} imply
$$
 \int_{E_j}\|\widetilde f(y,\cdot)\|_{L^2(\R_+)}\,d\mu_\kappa(y)
 \leq C \Xnorm{f} R_j^{\mathbf N}(1+\log R_j).
$$
For $p\geq0$, substitute $v=R_j^2/\tau$ to obtain
\[
 \int_0^T\tau^{-p}e^{-cR_j^2/\tau}\,d\tau
 = R_j^{2-2p}\int_{R_j^2/T}^\infty v^{p-2}e^{-cv}\,dv
 \leq C_p R_j^{2-2p}e^{-cR_j^2/(2T)},
 \qquad R_j^2/T\geq4.
\]
Minkowski's inequality, followed by this bound with $p=\ell+\mathbf N/2$, now implies
\begin{align}\label{e:tail2}
 &\sup_{\|x\|\leq R_0}
 \left\|\int_0^t\int_{\|y\|>R}|q_{t-s}(x,y)\widetilde f(y,s)|\,d\mu_\kappa(y)\,ds
 \right\|_{L^2(0,T)}\nonumber\\
 \leq&C \Xnorm{f}\sum_{j=0}^\infty R_j^{\mathbf N}(1+\log R_j)
 \int_0^T\tau^{-\ell-\mathbf N/2}e^{-cR_j^2/\tau}\,d\tau\nonumber\\
 \leq&C \Xnorm{f}\sum_{j=0}^\infty R_j^{2-2\ell}(1+\log R_j)e^{-cR_j^2/(2T)}\nonumber\\
 \leq&C_T \Xnorm{f} e^{-c'R^2/T}\to0.
\end{align}
In the last step, powers and logarithms are absorbed by half of the Gaussian,
and $R_j^2/T=4^jR^2/T\geq4^{j+1}$.
The same bound holds after restricting $t-s$ to any measurable subset of $(0,T)$.

For $q_\tau=h_\tau$, we use Cauchy--Schwarz in time on each annulus
and set $p=\mathbf N$ in the preceding integral bound:
\begin{align}\label{e:tail0}
 &\sup_{\|x\|\leq R_0}\sup_{0<t<T}
 \int_0^t\int_{\|y\|>R}h_{t-s}(x,y)|\widetilde f(y,s)|\,d\mu_\kappa(y)\,ds\nonumber\\
 \leq&C \Xnorm{f}\sum_{j=0}^\infty R_j^{\mathbf N}(1+\log R_j)
 \left(\int_0^T\tau^{-\mathbf N}e^{-cR_j^2/\tau}\,d\tau\right)^{1/2}\nonumber\\
 \leq&C \Xnorm{f}\sum_{j=0}^\infty R_j(1+\log R_j)e^{-cR_j^2/(4T)}
 \leq C_T \Xnorm{f} e^{-c'R^2/T}\to0.
\end{align}
The asserted spatial $L^2$ convergence follows on multiplying this uniform bound
by $\mu_\kappa(B(0,R_0))^{1/2}$.
The proof is complete.
\end{proof}

\begin{proof}[Proof of Theorem~\ref{t:main}]
Write $f=\widetilde f+a$ as in Section~\ref{s:avg}. For $R\geq1$, define
$$
 f_R=\widetilde f\mathbf1_{B(0,R)},\qquad
 u_R(t)=\int_0^t e^{(t-s)\Delta_\kappa}f_R(s)\,ds.
$$
The cutoff $f_R$ belongs to $L^2(\mu_\kappa\times\R_+)$ by \eqref{e:ng}.
Since the semigroup is strongly continuous and contractive,
$u_R\in C([0,T];L^2(\mu_\kappa))$ and
$$
 \|u_R(t)\|_{L^2(\mu_\kappa)}
 \leq\int_0^t\|f_R(s)\|_{L^2(\mu_\kappa)}\,ds
 \leq\sqrt t\,\|f_R\|_{L^2(\mu_\kappa\times(0,t))}.
$$
If $S>R\geq2\max\{1,R_0,\sqrt T\}$, then
$|f_S-f_R|\leq|\widetilde f|\mathbf1_{\{\|x\|>R\}}$.
Applying \eqref{e:tail0} to this difference, we obtain
$$
 \sup_{0\leq t\leq T}\|u_S(t)-u_R(t)\|_{L^2(B(0,R_0),\mu_\kappa)}
 \leq C_{R_0,T}\Xnorm{f} e^{-c'R^2/T}.
$$
Thus $(u_R)$ converges in $C([0,T];L^2(B(0,R_0),\mu_\kappa))$ to a function $\widetilde u$ with zero initial value.
Passing to the limit in the kernel integral by \eqref{e:tail0}, we find
\begin{equation}\label{e:sol}
 \widetilde u(x,t)=\int_0^t\int_{\R^N}h_{t-s}(x,y)\widetilde f(y,s)\,d\mu_\kappa(y)\,ds.
\end{equation}

For each fixed $R$, \eqref{e:mult} and dominated convergence imply
$$
 \mathcal M_\kappa^{a_0,b_0}f_R\to m_R
 \quad\text{in }L^2(\mu_\kappa\times\R_+),\qquad a_0\to0+,\quad b_0\to\infty,
$$
where $m_R=\int_0^t\Delta_\kappa e^{(t-s)\Delta_\kappa}f_R(s)\,ds$
is understood as this limit.
Approximation by functions smooth in time with values in the domain of $\Delta_\kappa$, followed by
\eqref{e:l2}, proves
$$
 \partial_tu_R=f_R+m_R,\qquad \Delta_\kappa u_R=m_R.
$$
For $S>R\geq2\max\{1,R_0,\sqrt T\}$, apply \eqref{e:tail2} to $f_S-f_R$
with the integral restricted to $a_0<t-s<b_0$. Passing to the $L^2$ limit
as $a_0\to0+$ and $b_0\to\infty$, we obtain
$$
 \|m_S-m_R\|_{L^2(B(0,R_0)\times(0,T),d\mu_\kappa\,dt)}
 \leq C_{R_0,T}\Xnorm{f} e^{-c'R^2/T}\to0.
$$
Thus $(m_R)$ converges in local space-time $L^2$, and we denote its limit by $m$.

For $\varphi\in C_c^\infty(\R^N\times[0,\infty))$, the spatial support of
$\Delta_\kappa\varphi$ lies in the finite union of reflections of the spatial support of $\varphi$.
The local convergence just proved therefore permits
\begin{align*}
 \int_0^\infty\int_{\R^N}\widetilde u(-\partial_t\varphi-\Delta_\kappa\varphi)\,d\mu_\kappa\,dt
 =&\lim_{R\to\infty}\int_0^\infty\int_{\R^N}
 u_R(-\partial_t\varphi-\Delta_\kappa\varphi)\,d\mu_\kappa\,dt\\
 =&\lim_{R\to\infty}\int_0^\infty\int_{\R^N}f_R\varphi\,d\mu_\kappa\,dt\\
 =&\int_0^\infty\int_{\R^N}\widetilde f\varphi\,d\mu_\kappa\,dt.
\end{align*}
The same local convergence in $\Delta_\kappa u_R=m_R$ proves
$$
 \Delta_\kappa \widetilde u=m,\qquad \partial_t\widetilde u=\widetilde f+m
 \quad\text{in weighted distributions}.
$$
Both right-hand sides belong to $L^2(B(0,R_0)\times(0,T),d\mu_\kappa\,dt)$.

To estimate $m$ by Proposition~\ref{p:mr}, we will prove that, in
$L^2(B(0,R_0)\times(0,T),d\mu_\kappa\,dt)$,
\begin{equation}\label{e:lim}
 \lim_{\substack{a_0\to0+\,,\ b_0\to\infty}}
 \mathcal M_\kappa^{a_0,b_0}\widetilde f=m=\Delta_\kappa \widetilde u.
\end{equation}
Decompose $\widetilde f=f_R+(\widetilde f-f_R)$.
In the next inequality all $L^2$ norms are taken on $B(0,R_0)\times(0,T)$:
\[
 \|\mathcal M_\kappa^{a_0,b_0}\widetilde f-m\|_2
 \leq \|\mathcal M_\kappa^{a_0,b_0}(\widetilde f-f_R)\|_2
 +\|\mathcal M_\kappa^{a_0,b_0}f_R-m_R\|_2+\|m_R-m\|_2.
\]
By \eqref{e:tail2} and the convergence $m_R\to m$, the first and last terms can be made arbitrarily small
by choosing $R$ sufficiently large, uniformly in $a_0,b_0$.
With $R$ fixed, the middle term tends to zero as $a_0\to0+$ and $b_0\to\infty$, proving \eqref{e:lim}.

For every ball $B$ and $T<\infty$, the map $g\mapsto g-(g(\cdot,t))_B$ is continuous on $L^2(B\times(0,T),d\mu_\kappa\,dt)$.
Hence \eqref{e:lim} and Proposition~\ref{p:mr} imply
\begin{align*}
 \frac1{\mu_\kappa(B)}\int_0^T\int_B
 |m(x,t)-(m(\cdot,t))_B|^2\,d\mu_\kappa(x)\,dt
 &\leq\liminf_{\substack{a_0\to0+\\b_0\to\infty}}
 \Xnorm{\mathcal M_\kappa^{a_0,b_0}\widetilde f}^2\\
 &\leq C\Xnorm{f}^2.
\end{align*}
Letting $T\to\infty$ and taking the supremum over $B$, we obtain
\[
 \Xnorm{\Delta_\kappa \widetilde u}= \Xnorm{m}\leq C\Xnorm{f},
\]
and
\[
 \Xnorm{\partial_t\widetilde u}\leq \Xnorm{\widetilde f}+\Xnorm{m}\leq C\Xnorm{f}.
\]

Since $a\in L^1(0,T)$ on every finite interval, define
$$
 u(x,t)=\widetilde u(x,t)+\int_0^t a(s)\,ds.
$$
Since $\int h_\tau(x,y)\,d\mu_\kappa(y)=1$, this is the Duhamel integral for $f$.
The added term is continuous, vanishes at $t=0$, and has zero Dunkl Laplacian. Thus
$$
 \Delta_\kappa u=m,\qquad \partial_tu=\widetilde f+a+m=f+m,
$$
and \eqref{e:mr} follows from the estimates for $\widetilde u$.

Let $v$ be another weak solution of \eqref{e:pde} in the class \eqref{e:class}, and write $w=v-u$.
Then $w(\cdot,0)=0$ and $\partial_tw=\Delta_\kappa w$ in weighted distributions for $t>0$.
By Lemma~\ref{l:rec}, for every $R,T>0$,
$$
 \partial_tw=\Delta_\kappa w
 \in L^2(B(0,R)\times(0,T),d\mu_\kappa\,dt),
 \qquad \Xnorm{\partial_tw}<\infty.
$$

The equation and the continuity of $w$ in time with values in $L^2_{\mathrm{loc}}(\mu_\kappa)$ imply that $w$ belongs to
$W^{1,2}(0,T;L^2(B(0,R),\mu_\kappa))$ for every $R,T>0$.
Hence $w(t)=\int_0^t\partial_sw(s)\,ds$ in local spatial $L^2$.
Write $b(t)=(\partial_tw(\cdot,t))_{B_*}$.
Since $b\in L^2(0,T)$, Lemma~\ref{l:norm}, applied to $\partial_tw$, implies
\begin{align}\label{e:ug}
 \sup_{0\le t\le T}\|w(t)\|_{L^2(B(0,R),\mu_\kappa)}
 \leq&\sqrt T\,
 \|\partial_tw-b\|_{L^2(B(0,R)\times\R_+,d\mu_\kappa\,dt)}
   +\mu_\kappa(B(0,R))^{1/2}\int_0^T|b(s)|\,ds\nonumber\\
 \leq&C_T R^{\mathbf N/2}(1+\log R),\qquad R\ge1.
\end{align}

The growth bound \eqref{e:ug} allows us to test the equation, after a spatial cutoff, against a solution of the backward heat equation.
Fix $t_0>0$, $\varphi\in C_c^\infty(\R^N)$, and $0<\varepsilon<1$.
For $0\le t\le t_0$, define
$$
 \psi_\varepsilon(x,t)
 =e^{(t_0-t+\varepsilon)\Delta_\kappa}\varphi(x).
$$
It satisfies $\partial_t\psi_\varepsilon+\Delta_\kappa\psi_\varepsilon=0$.
If $\operatorname{supp}\varphi\subset B(0,R_\varphi)$ and $\|x\|\ge2R_\varphi+2$,
Lemma~\ref{l:heat} and $d(x,y)\geq\|x\|-\|y\|$ imply
$$
 |\psi_\varepsilon(x,t)|+\|\nabla\psi_\varepsilon(x,t)\|
 \le C_{\varphi,t_0}e^{-c\|x\|^2/(t_0+1)},
 \qquad 0\le t\le t_0,\quad 0<\varepsilon<1.
$$
Indeed, the factors $\tau^{-\mathbf N/2}$ and
$\tau^{-(\mathbf N+1)/2}$ in the heat estimates are absorbed into half of the Gaussian,
uniformly for $0<\tau\le t_0+1$ and $\|x\|\ge1$.

Choose a radial $\chi\in C_c^\infty(\R^N)$ which is one on $B(0,1)$
and zero outside $B(0,2)$, and write $\chi_R(x)=\chi(x/R)$.
Then
$$
 \|\nabla\chi_R\|\le C R^{-1},\qquad
 |\Delta_\kappa\chi_R|\le C R^{-2}.
$$
Both derivatives are supported in $B(0,2R)\setminus B(0,R)$.
The second bound follows from the radial formula for $\Delta_\kappa$.
Since $\chi_R$ is $G$-invariant, the product formula has the form
\begin{equation}\label{e:prod}
 \Delta_\kappa(\chi_R\psi_\varepsilon)
 =\chi_R\Delta_\kappa\psi_\varepsilon
   +\psi_\varepsilon\Delta_\kappa\chi_R
   +2\nabla\chi_R\cdot\nabla\psi_\varepsilon.
\end{equation}
Testing the weak equation with $\chi_R\psi_\varepsilon$ and smooth time cutoffs on $[0,t_0]$,
then using local $L^2$ continuity and $w(0)=0$, we obtain
$$
 \int_{\R^N}w(x,t_0)\chi_R(x)e^{\varepsilon\Delta_\kappa}\varphi(x)\,d\mu_\kappa(x)
 =\int_0^{t_0}\int_{\R^N}w
   \bigl(\psi_\varepsilon\Delta_\kappa\chi_R
          +2\nabla\chi_R\cdot\nabla\psi_\varepsilon\bigr)\,d\mu_\kappa\,dt.
$$
Using \eqref{e:ug} for $w$, we estimate the right-hand side by
\begin{align*}
 &\left|\int_0^{t_0}\int_{\R^N}w
 \bigl(\psi_\varepsilon\Delta_\kappa\chi_R
 +2\nabla\chi_R\cdot\nabla\psi_\varepsilon\bigr)\,d\mu_\kappa\,dt\right|\\
 &\qquad\leq C_{\varphi,t_0}(R^{-1}+R^{-2})e^{-cR^2/(t_0+1)}
 \mu_\kappa(B(0,2R))^{1/2}
 \int_0^{t_0}\|w(t)\|_{L^2(B(0,2R),\mu_\kappa)}\,dt\\
 &\qquad\leq C_{\varphi,t_0}R^{\mathbf N}(1+\log R)e^{-cR^2/(t_0+1)}\to0,
\end{align*}
with convergence uniform for $0<\varepsilon<1$.
Also, $w(\cdot,t_0)e^{\varepsilon\Delta_\kappa}\varphi$ is integrable with respect to $\mu_\kappa$:
the Gaussian decay dominates the growth in \eqref{e:ug}.
We may therefore let $R\to\infty$ in the identity and conclude that
$$
 \int_{\R^N}w(x,t_0)e^{\varepsilon\Delta_\kappa}\varphi(x)\,d\mu_\kappa(x)=0.
$$
Take $R>2R_\varphi+2$. An application of Cauchy--Schwarz on each $B(0,2^{j+1}R)$ shows that
\begin{align*}
  \left|\int_{\R^N}w(x,t_0)
 \bigl(e^{\varepsilon\Delta_\kappa}\varphi(x)-\varphi(x)\bigr)\,d\mu_\kappa(x)\right|
 \leq&\|w(t_0)\|_{L^2(B(0,R),\mu_\kappa)}
 \|e^{\varepsilon\Delta_\kappa}\varphi-\varphi\|_{L^2(\mu_\kappa)}\\
 &+C_{\varphi,t_0}\sum_{j=0}^\infty
 (2^jR)^{\mathbf N}(1+\log(2^jR))e^{-c4^jR^2/(t_0+1)}.
\end{align*}
For fixed $R$, strong continuity on $L^2(\mu_\kappa)$ makes the first term tend to zero as
$\varepsilon\to0+$. The series has no dependence on $\varepsilon$ and tends to zero as $R\to\infty$.
Thus $\int w(x,t_0)\varphi(x)\,d\mu_\kappa(x)=0$ for every test function and every $t_0>0$,
which proves $w=0$.
This completes the proof of Theorem~\ref{t:main}.
\end{proof}

\section{Square functions and time-weighted energy}\label{sec:4}

\subsection{Centered square functions}\label{s:sq}

Lemma~\ref{l:heat} covers the ordinary gradient kernel. For the reflection quotient,
we integrate along the segment from $x$ to $\sigma_\alpha x$ near the hyperplane
and estimate the two kernel values separately farther away.

\begin{lemma}\label{l:ker}
For every $\alpha\in R_+$, $\tau>0$, $x,y\in\R^N$ and $|\beta|\leq1$,
\begin{equation}\label{e:ker}
 \left|\partial_x^\beta
 \frac{h_\tau(x,y)-h_\tau(\sigma_\alpha x,y)}{\langle\alpha,x\rangle}\right|
 \leq\frac{C\tau^{-(1+|\beta|)/2}e^{-c d(x,y)^2/\tau}}
 {\mu_\kappa(B(x,\sqrt\tau))(1+\|x-y\|/\sqrt\tau)}.
\end{equation}
The quotient extends smoothly across $\langle\alpha,x\rangle=0$.
The quotient and its first derivatives in $x$ have integral zero with respect to $d\mu_\kappa(y)$.
\end{lemma}

\begin{proof}
Fix $\alpha\in R_+$ and write $s=\sqrt\tau$, $p=\langle\alpha,x\rangle$,
$a=|p|/\|\alpha\|$, and $E=\|x-y\|$.

Suppose first that $a\leq s$.
For $z_\theta=x-2\theta p\alpha/\|\alpha\|^2$, the fundamental theorem of calculus takes the form
\begin{equation}\label{e:path}
 \frac{h_\tau(x,y)-h_\tau(\sigma_\alpha x,y)}{\langle\alpha,x\rangle}
 =\frac2{\|\alpha\|^2}\int_0^1\langle\alpha,\nabla_xh_\tau(z_\theta,y)\rangle\,d\theta.
\end{equation}
Here $\|z_\theta-x\|\leq2s$, so
$$
 \mu_\kappa(B(z_\theta,s))\simeq\mu_\kappa(B(x,s)),\qquad
 1+\frac{\|z_\theta-y\|}{s}\simeq1+\frac Es,\qquad
 d(z_\theta,y)\geq d(x,y)-2s.
$$
Moreover, the differential of $x\mapsto z_\theta$ is
$I-2\theta\alpha\otimes\alpha/\|\alpha\|^2$, whose operator norm is at most one.
For $|\beta|=0,1$, Lemma~\ref{l:heat} and \eqref{e:path} therefore imply
\[
 \left|\partial_x^\beta
 \frac{h_\tau(x,y)-h_\tau(\sigma_\alpha x,y)}{\langle\alpha,x\rangle}\right|
 \leq C_\alpha\int_0^1\|\nabla_x^{|\beta|+1}h_\tau(z_\theta,y)\|\,d\theta
 \leq \frac{C_\alpha s^{-1-|\beta|}e^{-c d(x,y)^2/s^2}}
 {\mu_\kappa(B(x,s))(1+E/s)^2}.
\]
Formula \eqref{e:path} also defines the quotient and its derivatives when $p=0$.

Suppose now that $a>s$, and write $E'=\|\sigma_\alpha x-y\|$.
Reflection invariance implies
$$
 \mu_\kappa(B(\sigma_\alpha x,s))=\mu_\kappa(B(x,s)),\qquad
 d(\sigma_\alpha x,y)=d(x,y),\qquad |E-E'|\leq2a.
$$
If $E\leq4a$, then $s+E\leq5a$ and
$$
 \frac1a\left\{(1+E/s)^{-2}+(1+E'/s)^{-2}\right\}
 \leq\frac2a\leq\frac{10}{s+E}.
$$
If $E>4a$, then $E'\geq E/2$ and
$$
 \frac1a\left\{(1+E/s)^{-2}+(1+E'/s)^{-2}\right\}
 \leq\frac Cs(1+E/s)^{-2}\leq\frac C{s+E}.
$$
Thus, in both cases,
\begin{equation}\label{e:end}
 \frac1a\left\{(1+E/s)^{-2}+(1+E'/s)^{-2}\right\}
 \leq\frac Cs(1+E/s)^{-1}.
\end{equation}
Applying the heat-kernel bound at $x$ and $\sigma_\alpha x$, we obtain
\[
 \left|\frac{h_\tau(x,y)-h_\tau(\sigma_\alpha x,y)}{\langle\alpha,x\rangle}\right|
 \leq \frac{C_\alpha e^{-c d(x,y)^2/s^2}}{\mu_\kappa(B(x,s))}
 \frac{(1+E/s)^{-2}+(1+E'/s)^{-2}}a
 \leq \frac{C_\alpha s^{-1}e^{-c d(x,y)^2/s^2}}{\mu_\kappa(B(x,s))(1+E/s)}.
\]
For the derivative, the quotient rule is
\[
 \nabla_x\frac{h_\tau(x,y)-h_\tau(\sigma_\alpha x,y)}{\langle\alpha,x\rangle}
 = \frac{\nabla_xh_\tau(x,y)-\sigma_\alpha\nabla_xh_\tau(\sigma_\alpha x,y)}
 {\langle\alpha,x\rangle}
  -\frac{\alpha\{h_\tau(x,y)-h_\tau(\sigma_\alpha x,y)\}}{\langle\alpha,x\rangle^2}.
\]
Consequently, by \eqref{e:end} and $a^{-1}\leq s^{-1}$,
\begin{align*}
  \left\|\nabla_x\frac{h_\tau(x,y)-h_\tau(\sigma_\alpha x,y)}{\langle\alpha,x\rangle}\right\|
 \leq&\frac{C_\alpha e^{-c d(x,y)^2/s^2}}{\mu_\kappa(B(x,s))}
 \left(\frac1{as}+\frac1{a^2}\right)
 \left\{(1+E/s)^{-2}+(1+E'/s)^{-2}\right\}\\
 \leq&\frac{C_\alpha s^{-2}e^{-c d(x,y)^2/s^2}}{\mu_\kappa(B(x,s))(1+E/s)}.
\end{align*}
Taking the reflection quotient of $\int h_\tau(x,y)\,d\mu_\kappa(y)=1$
and differentiating in $x$ proves both cancellation identities. The bounds above
allow differentiation under the integral also at $p=0$.
The proof is complete.
\end{proof}

\begin{remark}\label{r:sharp}
The Euclidean exponent $1$ in \eqref{e:ker} is optimal already for $\beta=0$.
Let $R=\{-1,1\}$, $\kappa>0$, and write $s=\sqrt\tau$ for $\tau>0$.
For $a>0$, the rank-one estimates in \cite[Proposition~2.3(e) and Remark~2.4]{ABDH} imply
$$
 h_\tau(-a,-a)\simeq\frac1{\mu_\kappa(B(a,s))},\qquad
 h_\tau(a,-a)\simeq\frac1{\mu_\kappa(B(a,s))}\frac{s^2}{s^2+a^2}.
$$
For $a\geq Ms$, with $M$ sufficiently large depending only on $\kappa$,
$h_\tau(a,-a)\leq h_\tau(-a,-a)/2$. Hence
$$
 \left|\frac{h_\tau(a,-a)-h_\tau(-a,-a)}{a}\right|
 \simeq\frac1{a\mu_\kappa(B(a,s))}
 \simeq\frac{s^{-1}}{\mu_\kappa(B(a,s))(1+2a/s)}.
$$
An estimate of the form \eqref{e:ker} with Euclidean exponent $\eta>1$
would therefore imply $1\lesssim(s/a)^{\eta-1}$, which fails as $a/s\to\infty$.
Thus the reflection quotient cannot share the Euclidean exponent $2$ available
for the ordinary gradient in Lemma~\ref{l:heat}.
\end{remark}

\begin{proof}[Proof of Theorem~\ref{t:sq}]
We first prove the upper estimates.
For a Hilbert-valued function $v$, write
$$
 \Gamma_\kappa(v)
 =\|\nabla v\|_{\mathcal H^N}^2
 +\sum_{\alpha\in R_+}\frac{\kappa(\alpha)\|\alpha\|^2}{2}
 \left\|\frac{v-v\circ\sigma_\alpha}{\langle\alpha,\cdot\rangle}\right\|_{\mathcal H}^2.
$$
Apply \eqref{e:form} to finite orthonormal expansions in $\mathcal H$ and pass to the closed form.
For $V\in L^2(\mu_\kappa;\mathcal H)$ and $0<a<b<\infty$,
\begin{align*}
 \int_a^b\int_{\R^N}\Gamma_\kappa(e^{t\Delta_\kappa}V)\,d\mu_\kappa\,dt
 &=-\operatorname{Re}\int_a^b
 \langle\Delta_\kappa e^{t\Delta_\kappa}V,e^{t\Delta_\kappa}V\rangle_{L^2(\mu_\kappa;\mathcal H)}\,dt\\
 &=\frac12\left(\|e^{a\Delta_\kappa}V\|_2^2-\|e^{b\Delta_\kappa}V\|_2^2\right)
 \leq\frac12\|V\|_2^2.
\end{align*}
Let $a\to0+$ and $b\to\infty$. Then
\begin{equation}\label{e:el2}
 \int_0^\infty\int_{\R^N}\Gamma_\kappa(e^{t\Delta_\kappa}V)\,d\mu_\kappa\,dt
 \leq\frac12\|V\|_2^2.
\end{equation}

For $0<a<b<\infty$, let $K_t(x,y)$ denote either of the kernels
$$
 \nabla_xh_t(x,y),\qquad
 \sqrt{\kappa(\alpha)/2}\|\alpha\|
 \frac{h_t(x,y)-h_t(\sigma_\alpha x,y)}{\langle\alpha,x\rangle}
 \quad(\alpha\in R_+,\ \kappa(\alpha)>0).
$$
For bounded compactly supported $F$, put
$\mathcal T_tF(x)=\int K_t(x,y)F(y)\,d\mu_\kappa(y)$.
The two choices produce the ordinary heat gradient and the weighted reflection quotient.
Estimate \eqref{e:el2} bounds the $L^2$ extension of $F\mapsto(\mathcal T_tF)_{a<t<b}$ by $1/\sqrt2$,
with values in $L^2((a,b);\mathcal H^N)$ and $L^2((a,b);\mathcal H)$, respectively.
Both operators vanish on constants by the conservation property and Lemma~\ref{l:ker}.
For $|\beta|\leq1$, Lemmas~\ref{l:heat} and \ref{l:ker} imply
$$
 \|\partial_x^\beta K_t(x,y)\|
 \leq\frac{Ct^{-(1+|\beta|)/2}e^{-cd(x,y)^2/t}}
 {\mu_\kappa(B(x,\sqrt t))(1+\|x-y\|/\sqrt t)}.
$$
Lemma~\ref{l:int}, with $p=2$ and $\ell=|\beta|$, bounds the time integral by
$$
 \left(\int_a^b\|\partial_x^\beta K_t(x,y)\|^2\,dt\right)^{1/2}
 \leq\frac{Cd(x,y)^{-|\beta|}}{\mu_\kappa(B(x,d(x,y)))}
 \left(\frac{d(x,y)}{\|x-y\|}\right)^{1/2}.
$$
By \eqref{e:kd}, these bounds imply the kernel conditions in Proposition~\ref{p:cz}.
Lemma~\ref{l:ker} allows the segment in \eqref{e:kd} to cross a reflecting hyperplane.
Keeping $a,b$ fixed in the pointwise bound, we also have
$$
 \left(\int_a^b\|K_t(x,y)\|^2\,dt\right)^{1/2}
 \leq C_{a,b}e^{-c d(x,y)^2/b}.
$$
For $F\in\mathrm{BMO}_2(\mu_\kappa;\mathcal H)$, Lemma~\ref{l:conv} ensures absolute convergence of the spatial integral
in the $L^2$ norm over $a<t<b$. If $F$ also belongs to $L^2(\mu_\kappa;\mathcal H)$,
the integral represents the $L^2$ operator.
Proposition~\ref{p:cz} therefore proves
$$
 \sup_B\left(\frac1{\mu_\kappa(B)}\int_B\int_a^b
 \|\mathcal T_tF(x)-(\mathcal T_tF)_B\|^2
 \,dt\,d\mu_\kappa(x)\right)^{1/2}
 \leq C\Bnorm{F},
$$
uniformly in $a,b$.
By monotone convergence, we obtain \eqref{e:sq} and \eqref{e:rsq}.
The upper bound in \eqref{e:dsq} then follows from the definition of $T_j$.

For the reverse estimate, we will use the adjoint of the truncated Dunkl gradient to obtain the difference $e^{2a\Delta_\kappa}F-e^{2b\Delta_\kappa}F$.
Fix $0<a<b<\infty$ and write
$$
 \mathcal H_{a,b}=L^2((a,b),dt;\mathcal H^N),\qquad
 S_{a,b}F(x)=\bigl(T_je^{t\Delta_\kappa}F(x)\bigr)_{j=1}^N,
 \quad a<t<b.
$$
The identity $T_j^*=-T_j$ determines the adjoint.
The energy calculation above, with $\Delta_\kappa=\sum_jT_j^2$, bounds its $L^2$ norm:
\begin{equation}\label{e:adj}
 S_{a,b}^*G=-\int_a^b\sum_{j=1}^N T_je^{t\Delta_\kappa}G_j(\cdot,t)\,dt,
 \qquad \|S_{a,b}^*\|_{2\to2}=\|S_{a,b}\|_{2\to2}\leq2^{-1/2}.
\end{equation}
By \cite[Lemma~3.3]{DH1},
$$
 T_{j,x}h_t(x,y)=\frac{y_j-x_j}{2t}h_t(x,y)=-T_{j,y}h_t(x,y).
$$
Differentiating this identity and applying Lemma~\ref{l:heat}, for $|\beta|+|\gamma|\leq1$ we obtain
\begin{equation}\label{e:dk}
 |\partial_x^\beta\partial_y^\gamma T_{j,x}h_t(x,y)|
 \leq\frac{Ct^{-(1+|\beta|+|\gamma|)/2}e^{-cd(x,y)^2/t}}
 {\mu_\kappa(B(x,\sqrt t))(1+\|x-y\|/\sqrt t)}.
\end{equation}
The kernel of \eqref{e:adj} is
$$
 K_{a,b}^*(x,y)g=-\int_a^b\sum_jT_{j,x}h_t(x,y)g_j(t)\,dt,
 \qquad g\in\mathcal H_{a,b}.
$$
For $d=d(x,y)>0$, $E=\|x-y\|$ and $|\beta|\leq1$,
Cauchy--Schwarz, \eqref{e:dk} and Lemma~\ref{l:int} imply
\[
 \|\partial_x^\beta K_{a,b}^*(x,y)\|
 \leq \left(\int_a^b\sum_j|\partial_x^\beta T_{j,x}h_t(x,y)|^2\,dt\right)^{1/2}
 \leq \frac{Cd^{-|\beta|}}{\mu_\kappa(B(x,d))}\left(\frac dE\right)^{1/2}.
\]
Thus \eqref{e:kd} verifies the kernel conditions of Proposition~\ref{p:cz}, uniformly in $a,b$.
The conservation property also implies, for every spatially constant $g\in\mathcal H_{a,b}$,
$$
 \int_{\R^N}K_{a,b}^*(x,y)g\,d\mu_\kappa(y)
 =-\int_a^b\sum_jg_j(t)T_j1\,dt=0.
$$
For fixed $a,b$, \eqref{e:dk} implies
$$
 \|K_{a,b}^*(x,y)\|
 \leq\left(\int_a^b\sum_j|T_{j,x}h_t(x,y)|^2\,dt\right)^{1/2}
 \leq C_{a,b}e^{-c d(x,y)^2/b}.
$$
For bounded compactly supported $G:\R^N\to\mathcal H_{a,b}$, Fubini's theorem
identifies the integral with $S_{a,b}^*G$. By the Gaussian bound and Lemma~\ref{l:conv},
it converges absolutely for $G\in\mathrm{BMO}_2(\mu_\kappa;\mathcal H_{a,b})$ and agrees with $S_{a,b}^*G$ for every $G\in L^2(\mu_\kappa;\mathcal H_{a,b})$.
Proposition~\ref{p:cz} and \eqref{e:adj} now imply
\begin{equation}\label{e:abmo}
 \Bnorm{S_{a,b}^*G}\leq C\|G\|_{\mathrm{BMO}_2(\mu_\kappa;\mathcal H_{a,b})},
\end{equation}
where $C$ has no dependence on $a,b$ or $\mathcal H$.
For $F_R=F\mathbf1_{B(0,R)}$, applying the Dunkl transform componentwise in $\mathcal H$, we obtain in $L^2$ that
$$
 S_{a,b}^*S_{a,b}F_R
 =-\int_a^b\Delta_\kappa e^{2t\Delta_\kappa}F_R\,dt
 =\frac12\left(e^{2a\Delta_\kappa}F_R-e^{2b\Delta_\kappa}F_R\right).
$$
For $\|x\|\leq R_0$ and $a\leq t\leq b$, \eqref{e:dk} implies
$$
 |T_{j,x}h_t(x,y)T_{j,y}h_t(y,z)|
 \leq C_{a,b,R_0}e^{-c_{a,b}(\|y\|^2+\|z\|^2)}.
$$
Indeed, $\|y\|\leq\|x\|+d(x,y)$ and $\|z\|\leq\|y\|+d(y,z)$.
The growth bound used in Lemma~\ref{l:conv} therefore implies
\begin{align*}
 &\sup_{\|x\|\leq R_0}\int_a^b\int_{\R^N}\int_{\R^N}
 |T_{j,x}h_t(x,y)T_{j,y}h_t(y,z)|\|F(z)\|_{\mathcal H}
 \,d\mu_\kappa(z)\,d\mu_\kappa(y)\,dt\\
 &\qquad\leq C_{a,b,R_0}\int_{\R^N}
 e^{-c_{a,b}\|z\|^2}\|F(z)\|_{\mathcal H}\,d\mu_\kappa(z)<\infty.
\end{align*}
Fubini and dominated convergence imply, locally uniformly in $x$,
$$
 S_{a,b}^*S_{a,b}F_R\to S_{a,b}^*S_{a,b}F.
$$
Lemma~\ref{l:conv} proves the same convergence for the two heat terms.
Letting $R\to\infty$, we obtain
\begin{equation}\label{e:rep}
 S_{a,b}^*S_{a,b}F
 =\frac12\left(e^{2a\Delta_\kappa}F-e^{2b\Delta_\kappa}F\right).
\end{equation}
Combining \eqref{e:rep} with \eqref{e:abmo}, we obtain
\begin{equation}\label{e:rbmo}
 \Bnorm{e^{2a\Delta_\kappa}F-e^{2b\Delta_\kappa}F}
 \leq C\|S_{a,b}F\|_{\mathrm{BMO}_2(\mu_\kappa;\mathcal H_{a,b})}.
\end{equation}

As $a\to0+$, $e^{2a\Delta_\kappa}F\to F$ in local $L^2(\mu_\kappa;\mathcal H)$,
by strong continuity on spatial cutoffs and Gaussian control of the tails.
For the other endpoint, apply Lemma~\ref{l:mean} at the center $x$,
with both radii equal to $\sqrt t$ and $\eta=2$. By cancellation and Lemma~\ref{l:heat},
\[
 \|\nabla e^{t\Delta_\kappa}F(x)\|_{\mathcal H^N}
 = \left\|\int_{\R^N}\nabla_xh_t(x,y)
 \bigl(F(y)-F_{B(x,4\sqrt t)}\bigr)\,d\mu_\kappa(y)\right\|_{\mathcal H^N}
 \leq Ct^{-1/2}\Bnorm{F}.
\]
Hence, for every fixed ball $B=B(c,r)$,
\begin{align*}
 \left(\frac1{\mu_\kappa(B)}\int_B
 \|e^{2b\Delta_\kappa}F-(e^{2b\Delta_\kappa}F)_B\|_{\mathcal H}^2\,d\mu_\kappa\right)^{1/2}
 &\leq2r\|\nabla e^{2b\Delta_\kappa}F\|_{L^\infty(\mu_\kappa;\mathcal H^N)}\\
 &\leq Crb^{-1/2}\Bnorm{F}\to0\qquad(b\to\infty).
\end{align*}
Take these limits on each fixed ball in \eqref{e:rbmo}, and then take the supremum over balls.
By monotone convergence in $a,b$, we obtain
$$
 \Bnorm{F}
 \leq C\sup_{0<a<b<\infty}\|S_{a,b}F\|_{\mathrm{BMO}_2(\mu_\kappa;\mathcal H_{a,b})}
 =C\Xnorm{\nabla_\kappa e^{t\Delta_\kappa}F}.
$$
This completes the proof of Theorem~\ref{t:sq}.
\end{proof}

\begin{remark}\label{r:ctr}
The subtraction of the spatial mean in \eqref{e:dsq} is necessary at large times.
Even for $\kappa=0$ and $v(x)=\operatorname{sgn}x$ on $\R$,
$$
 \partial_xe^{\tau\Delta}v(x)
 =(\pi\tau)^{-1/2}e^{-x^2/(4\tau)},
 \qquad
 \int_1^\infty|\partial_xe^{\tau\Delta}v(x)|^2\,d\tau=\infty.
$$

The centered heat-gradient seminorm can vanish for a function outside BMO.
For $F(x)=x_\ell$, $1\leq\ell\leq N$,
$$
 e^{t\Delta_\kappa}F=F,\qquad
 T_jF=\delta_{j\ell}+2\sum_{\alpha\in R_+}
 \kappa(\alpha)\frac{\alpha_j\alpha_\ell}{\|\alpha\|^2}.
$$
Thus $\Xnorm{\nabla_\kappa e^{t\Delta_\kappa}F}=0$, although $F\notin\mathrm{BMO}_2(\mu_\kappa)$.
\end{remark}

\subsection{Time-weighted energy}\label{s:en}

For space-time functions, reflection acts on the spatial variable:
$(u\circ\sigma_\alpha)(x,t)=u(\sigma_\alpha x,t)$.

\begin{corollary}\label{c:en}
For $f$ and $u$ as in Theorem~\ref{t:main},
\begin{equation}\label{e:en}
 \Xnorm{t^{-1/2}\nabla u}\leq C\Xnorm{f}.
\end{equation}
For each $\alpha\in R_+$ with $\kappa(\alpha)>0$,
\begin{equation}\label{e:ren}
 \biggl\|t^{-1/2}\frac{u-u\circ\sigma_\alpha}{\langle\alpha,\cdot\rangle}\biggr\|_{\mathbb L^2(\R_+;\mathrm{BMO}(\mu_\kappa))}
 \leq C\Xnorm{f}.
\end{equation}
\end{corollary}

\begin{proof}
Let $\widetilde f$ and $\widetilde u$ be as in Lemma~\ref{l:norm} and \eqref{e:sol},
and regard $\widetilde f$ as a function of $x$ with values in $\mathcal H=L^2(\R_+,ds)$.
Fix $R_0,T>0$. Since $(\widetilde f(\cdot,s))_{B(0,1)}=0$, \eqref{e:av} implies
$$
 \|(\widetilde f(\cdot,s))_{B(x,4\sqrt\tau)}\|_{L^2(\R_+,ds)}
 \leq C_{R_0,T}(1+|\log\tau|)\Xnorm{f},
 \qquad \|x\|\leq R_0,\quad 0<\tau\leq T.
$$
Decompose $\widetilde f$ around this mean and apply Lemma~\ref{l:mean}
at the center $x$, with both radii equal to $\sqrt\tau$ and $\eta=2$.
Together with \eqref{e:heat} and the Gaussian bound
$\int|\partial_x^\nu h_\tau(x,y)|\,d\mu_\kappa(y)\leq C\tau^{-|\nu|/2}$, this proves, for $|\nu|\leq1$,
\begin{equation}\label{e:nh}
\begin{split}
 \int_{\R^N}|\partial_x^\nu h_\tau(x,y)|\|\widetilde f(y,\cdot)\|_{\mathcal H}\,d\mu_\kappa(y)
 &\leq C\tau^{-|\nu|/2}\bigl(\Xnorm{f}+\|(\widetilde f(\cdot,s))_{B(x,4\sqrt\tau)}\|_{L^2(\R_+,ds)}\bigr)\\
 &\leq C_{R_0,T}\tau^{-|\nu|/2}(1+|\log\tau|)\Xnorm{f}.
\end{split}
\end{equation}

The bound \eqref{e:nh} is integrable in $\tau$ for $|\nu|\leq1$.
By \eqref{e:sol} and contractivity of time translation on $L^2(\R_+)$,
differentiation of the Bochner integral proves $\widetilde u\in C^1(\R^N;L^2(0,T))$, with
$$
 \nabla \widetilde u(x,t)
 =\int_0^t\nabla e^{(t-s)\Delta_\kappa}\widetilde f(\cdot,s)(x)\,ds
 \quad\text{in }L^2(0,T).
$$
For $n_\alpha=\alpha/\|\alpha\|$, the fundamental theorem of calculus in $L^2(0,T)$ implies
$$
 \frac{\widetilde u(x,\cdot)-\widetilde u(\sigma_\alpha x,\cdot)}{\langle\alpha,x\rangle}
 =\frac2{\|\alpha\|}\int_0^1
 \partial_{n_\alpha}\widetilde u(\sigma_\alpha x+\theta(x-\sigma_\alpha x),\cdot)\,d\theta.
$$
This formula also defines the quotient on the hyperplane. By \eqref{e:nh} and Fubini's theorem,
we may take the reflection quotient inside the Duhamel integral for $\widetilde u$.

Fix a ball $B$. By \eqref{e:sq}, the difference between
$\nabla e^{\tau\Delta_\kappa}\widetilde f(\cdot,s)$ and its spatial mean on $B$ belongs to
$L^2(B\times\R_+\times\R_+,d\mu_\kappa(x)\,d\tau\,ds)$.
Cauchy--Schwarz in the Duhamel integral, the substitution $\tau=t-s$, and Tonelli's theorem imply
\begin{align*}
 &\frac1{\mu_\kappa(B)}\int_B\int_0^\infty
 t^{-1}\|\nabla \widetilde u(x,t)-(\nabla \widetilde u(\cdot,t))_B\|^2\,dt\,d\mu_\kappa(x)\\
 \leq&\frac1{\mu_\kappa(B)}\int_B\int_0^\infty\int_0^t
 \left\|\nabla e^{(t-s)\Delta_\kappa}\widetilde f(\cdot,s)(x)
 -(\nabla e^{(t-s)\Delta_\kappa}\widetilde f(\cdot,s))_B\right\|^2\,ds\,dt\,d\mu_\kappa(x)\\
 =&\frac1{\mu_\kappa(B)}\int_B\int_0^\infty\int_0^\infty
 \left\|\nabla e^{\tau\Delta_\kappa}\widetilde f(\cdot,s)(x)
 -(\nabla e^{\tau\Delta_\kappa}\widetilde f(\cdot,s))_B\right\|^2\,d\tau\,ds\,d\mu_\kappa(x)\\
 \leq&C\Xnorm{f}^2.
\end{align*}
The same calculation, using \eqref{e:rsq}, applies to each signed reflection quotient with $\kappa(\alpha)>0$.
Taking the supremum over $B$ proves \eqref{e:en} and \eqref{e:ren} for $\widetilde u$.
Since $u-\widetilde u$ depends only on time, both estimates hold for $u$.
The proof is complete.
\end{proof}

The Dunkl gradient satisfies
$$
 \nabla_\kappa u=\nabla u+
 \sum_{\alpha\in R_+}\kappa(\alpha)\alpha
 \frac{u-u\circ\sigma_\alpha}{\langle\alpha,\cdot\rangle}.
$$
Corollary~\ref{c:en} therefore also implies $\Xnorm{t^{-1/2}\nabla_\kappa u}\leq C\Xnorm{f}$.

\section{Nonzero initial values}\label{sec:5}

To apply Theorem~\ref{t:sq} to $T_j u_0$, we first prove that $e^{t\Delta_\kappa}u_0$ has initial value $u_0$
in local $L^2$ and commutes with $T_j$.
For $g\in\mathrm{BMO}_2(\mu_\kappa)$, the concentric-ball comparison \eqref{e:av} implies
\begin{equation}\label{e:bg}
 \|g\|_{L^2(B(0,r),\mu_\kappa)}
 \leq Cr^{\mathbf N/2}
 \bigl(|g_{B(0,1)}|+(1+\log r)\|g\|_{\mathrm{BMO}_2(\mu_\kappa)}\bigr),\qquad r\geq1.
\end{equation}

\begin{lemma}\label{l:grow}
Let $u_0\in L^2_{\mathrm{loc}}(\mu_\kappa)$ have weak Dunkl derivatives
$T_j u_0\in\mathrm{BMO}_2(\mu_\kappa)$ for $1\leq j\leq N$. Then for some $C_{u_0}>0$,
$$
 \|u_0\|_{L^2(B(0,r),\mu_\kappa)}
 \leq C_{u_0}r^{\mathbf N/2+1}(1+\log r),\qquad r\geq1.
$$
\end{lemma}

\begin{proof}
On the unit ball $B_*$, write $v_r(x)=u_0(rx)$ and $g_j=T_j u_0$.
For $\phi\in C_c^\infty(B_*)$, the definition of $T_j$ implies
$$
 \sum_{j=1}^NT_j(x_j\phi)
 =x\cdot\nabla\phi+\mathbf N\phi
  -\sum_{\alpha\in R_+}\kappa(\alpha)(\phi-\phi\circ\sigma_\alpha).
$$
Change variables $y=rx$ in $\int_{B_*}v_r\phi\,d\mu_\kappa$ and use the weak identities $T_j u_0=g_j$.
We obtain, as an $L^2(B_*,\mu_\kappa)$-valued distribution in $r$,
$$
 \partial_r v_r+\frac1r\sum_{\alpha\in R_+}\kappa(\alpha)(v_r-v_r\circ\sigma_\alpha)
 =\sum_{j=1}^Nx_jg_j(r\,\cdot).
$$
Both terms other than $\partial_rv_r$ belong to $L^2((a,b);L^2(B_*,\mu_\kappa))$ whenever $0<a<b<\infty$.
Hence $r\mapsto v_r$ belongs to $W^{1,2}_{\mathrm{loc}}((0,\infty);L^2(B_*,\mu_\kappa))$.
Density and continuity of dilations show that $r\mapsto u_0(r\,\cdot)$ is continuous in $L^2(B_*,\mu_\kappa)$.
Thus its continuous Sobolev representative equals $u_0(r\,\cdot)$ for every $r>0$.
Taking the inner product with $v_r$ and using reflection invariance, we obtain
\begin{align*}
 \frac12\frac{d}{dr}\|v_r\|_{L^2(B_*,\mu_\kappa)}^2
 +\frac1{2r}\sum_{\alpha\in R_+}\kappa(\alpha)
 \|v_r-v_r\circ\sigma_\alpha\|_{L^2(B_*,\mu_\kappa)}^2
 =&\operatorname{Re}\int_{B_*}\sum_{j=1}^Nx_jg_j(rx)\overline{v_r(x)}\,d\mu_\kappa(x)\\
 \leq&r^{-\mathbf N/2}\sum_{j=1}^N\|g_j\|_{L^2(B(0,r),\mu_\kappa)}
 \|v_r\|_{L^2(B_*,\mu_\kappa)}.
\end{align*}
Since $\kappa\geq0$, integration from $1$ to $r\geq1$ and \eqref{e:bg} imply
\begin{align*}
 r^{-\mathbf N/2}\|u_0\|_{L^2(B(0,r),\mu_\kappa)}
 =&\|v_r\|_{L^2(B_*,\mu_\kappa)}\\
 \leq&\|u_0\|_{L^2(B_*,\mu_\kappa)}
 +\int_1^r\rho^{-\mathbf N/2}\sum_{j=1}^N
 \|g_j\|_{L^2(B(0,\rho),\mu_\kappa)}\,d\rho\\
 \leq&\|u_0\|_{L^2(B_*,\mu_\kappa)}
 +Cr\sum_{j=1}^N\bigl(|(g_j)_{B_*}|+(1+\log r)\|g_j\|_{\mathrm{BMO}_2(\mu_\kappa)}\bigr).
\end{align*}
The integration is justified first for $(\|v_r\|_2^2+\varepsilon^2)^{1/2}$ and then by letting
$\varepsilon\to0+$.
The proof is complete.
\end{proof}

\begin{lemma}\label{l:init}
Suppose $u_0\in L^2_{\mathrm{loc}}(\mu_\kappa)$ satisfies, for some $C_0,m\geq0$,
\begin{equation}\label{e:grow}
 \|u_0\|_{L^2(B(0,R_0),\mu_\kappa)}\leq C_0R_0^m,\qquad R_0\geq1.
\end{equation}
Then the integral
$$
 e^{t\Delta_\kappa}u_0(x)=\int_{\R^N}h_t(x,y)u_0(y)\,d\mu_\kappa(y)
$$
converges absolutely and defines a smooth solution of the homogeneous Dunkl heat equation for $t>0$.
This solution extends to $C([0,\infty);L^2_{\mathrm{loc}}(\mu_\kappa))$ with value $u_0$ at $t=0$.
For each $j$ such that the weak derivative $T_j u_0$ belongs to $\mathrm{BMO}_2(\mu_\kappa)$,
\begin{equation}\label{e:comm}
 T_j e^{t\Delta_\kappa}u_0=e^{t\Delta_\kappa}(T_j u_0),\qquad t>0.
\end{equation}
\end{lemma}

\begin{proof}
Take $R_0,T>0$ and $R>2\max\{1,R_0,\sqrt T\}$.
On each annulus, Cauchy--Schwarz and \eqref{e:grow}, with the homogeneity of $\mu_\kappa$, imply
$$
 \int_{2^\ell R<\|y\|<2^{\ell+1}R}|u_0(y)|\,d\mu_\kappa(y)
 \leq C_mC_0(2^\ell R)^{m+\mathbf N/2},\qquad \ell\geq0.
$$
On $\|x\|\leq R_0$, $\|y\|>R$, the orbit distance satisfies $d(x,y)\geq\|y\|/2$.
Insert this and $\mu_\kappa(B(x,\sqrt t))\geq ct^{\mathbf N/2}$ into Lemma~\ref{l:heat}.
The tail is then bounded by
\begin{equation}\label{e:itail}
 \sup_{0<t\leq T}\sup_{\|x\|\leq R_0}
 \int_{\|y\|>R}h_t(x,y)|u_0(y)|\,d\mu_\kappa(y)
 \leq C_mC_0\sum_{\ell=0}^\infty
 (2^\ell R)^{m-\mathbf N/2}e^{-c4^\ell R^2/T}
 \leq C_{m,T}C_0e^{-c'R^2/T}\to0.
\end{equation}
As in the proof of Lemma~\ref{l:tail}, the powers of $2^\ell R$ are absorbed by half of the Gaussian factor.
The same calculation for the derivatives of $h_t$ is uniform for $t$ in compact subintervals of $(0,\infty)$,
so we may differentiate under the integral to obtain the heat equation.
The function $e^{t\Delta_\kappa}(u_0\mathbf1_{B(0,R)})$ is continuous in $L^2(\mu_\kappa)$ for $t\geq0$.
By \eqref{e:itail}, these functions converge in
$C([0,T];L^2(B(0,R_0),\mu_\kappa))$ to a function with initial value $u_0$ on $B(0,R_0)$.

Suppose now that $T_j u_0\in\mathrm{BMO}_2(\mu_\kappa)$.
By \eqref{e:bg}, the heat-kernel integral defining $e^{t\Delta_\kappa}(T_ju_0)$ also converges absolutely.
By symmetry and Lemma~3.3 of \cite{DH1},
$$
 T_{j,x}h_t(x,y)=\frac{y_j-x_j}{2t}h_t(x,y)=-T_{j,y}h_t(x,y).
$$
Take the radial cutoff $\chi_R$ used in the proof of Theorem~\ref{t:main}.
Its $G$-invariance implies
$T_j(\chi_R\varphi)=\chi_RT_j\varphi+(\partial_j\chi_R)\varphi$.
The weak derivative identity, with test function $\chi_R(y)h_t(x,y)$, becomes
\[
 \int_{\R^N}\chi_R(y)h_t(x,y)T_j u_0(y)\,d\mu_\kappa(y)
 = \int_{\R^N}\chi_R(y)T_{j,x}h_t(x,y)u_0(y)\,d\mu_\kappa(y)
 -\int_{\R^N}(\partial_j\chi_R)(y)h_t(x,y)u_0(y)\,d\mu_\kappa(y).
\]
For fixed $t>0$ and $\|x\|\leq R_0$, the absolute value of the last integral is bounded by
$$
 C_{m,t,R_0}C_0R^{m+\mathbf N/2-1}e^{-cR^2/t}\to0.
$$
The other integrals converge absolutely by the growth bounds and the Gaussian estimates.
Applying $T_j$ to the heat-kernel integral now proves \eqref{e:comm}.
The proof is complete.
\end{proof}

\begin{proof}[Proof of Corollary~\ref{c:init}]
Lemma~\ref{l:grow} verifies the growth hypothesis of Lemma~\ref{l:init}.

Using \eqref{e:comm} and $\Delta_\kappa=\sum_jT_j^2$, we have
$$
 \partial_t e^{t\Delta_\kappa}u_0
 =\Delta_\kappa e^{t\Delta_\kappa}u_0
 =\sum_{j=1}^NT_j e^{t\Delta_\kappa}(T_j u_0),
$$
and the upper estimate in \eqref{e:dsq} implies
\begin{equation}\label{e:ih}
\begin{split}
 \Xnorm{\partial_t e^{t\Delta_\kappa}u_0}
 &=\Xnorm{\Delta_\kappa e^{t\Delta_\kappa}u_0}\\
 &\leq\sum_{j=1}^N\Xnorm{T_j e^{t\Delta_\kappa}(T_j u_0)}
 \leq C\sum_{j=1}^N\|T_j u_0\|_{\mathrm{BMO}_2(\mu_\kappa)}.
\end{split}
\end{equation}
By the local $L^2$ continuity at $t=0$ in Lemma~\ref{l:init},
we may apply Lemma~\ref{l:rec} with \eqref{e:ih} to obtain
$$
 \Delta_\kappa e^{t\Delta_\kappa}u_0
 \in L^2(B\times(0,T),d\mu_\kappa\,dt)
 \quad\text{for every ball }B\text{ and every }T>0.
$$

Let $u$ be the sum of $e^{t\Delta_\kappa}u_0$ and the solution of Theorem~\ref{t:main}.
Then $u$ is continuous in local spatial $L^2$ with initial value $u_0$.
Since $f(x,t)-(f(\cdot,t))_{B(0,1)}$ and $\Delta_\kappa u$ are locally square-integrable and
$(f(\cdot,t))_{B(0,1)}\in L^1(0,T)$, the equation and the continuity of $u$ in time with values in $L^2_{\mathrm{loc}}(\mu_\kappa)$ imply
$$
 u\in W^{1,1}(0,T;L^2(B(0,R_0),\mu_\kappa)).
$$
Test the equation against $\varphi\in C_c^\infty(\R^N\times[0,\infty))$ and integrate over $t>\varepsilon$:
\begin{align*}
 \int_\varepsilon^\infty\int_{\R^N}
 u(-\partial_t\varphi-\Delta_\kappa\varphi)\,d\mu_\kappa\,dt
 =&\int_\varepsilon^\infty\int_{\R^N}f\varphi\,d\mu_\kappa\,dt
 +\int_{\R^N}u(x,\varepsilon)\varphi(x,\varepsilon)\,d\mu_\kappa(x)\\
 \to&\int_0^\infty\int_{\R^N}f\varphi\,d\mu_\kappa\,dt
 +\int_{\R^N}u_0(x)\varphi(x,0)\,d\mu_\kappa(x)
 \qquad(\varepsilon\to0+).
\end{align*}
Since $u$ has local $L^2$ trace $u_0$ at $t=0$, the boundary term converges, and \eqref{e:weak} holds.
Combining \eqref{e:mr} and \eqref{e:ih} proves \eqref{e:imr}.
If $v$ is another solution with initial value $u_0$ in \eqref{e:class},
then $v-e^{t\Delta_\kappa}u_0$ solves \eqref{e:pde} in the same class.
Theorem~\ref{t:main} implies $v=u$.
This completes the proof of Corollary~\ref{c:init}.
\end{proof}

\begin{remark}\label{r:null}
For the initial value $u_0(x)=\|x\|^2$, we have
$$
 T_j u_0=2x_j,\qquad \Delta_\kappa u_0=2\mathbf N,\qquad
 e^{t\Delta_\kappa}u_0=\|x\|^2+2\mathbf Nt.
$$
Thus the homogeneous maximal regularity seminorm vanishes. The mean-square oscillation of $T_j u_0$ on $B(0,r)$
is a positive multiple of $r^2$, so the derivative assumption of Corollary~\ref{c:init} fails.

Whenever $u_0\in L^2_{\mathrm{loc}}(\mu_\kappa)$ satisfies \eqref{e:grow}, the identity
$$
 \Xnorm{\Delta_\kappa e^{t\Delta_\kappa}u_0}=0
$$
is equivalent to $\Delta_\kappa u_0$ being a spatial constant as a weighted distribution.
Indeed, a vanishing seminorm and smoothness at positive times imply
$\Delta_\kappa e^{t\Delta_\kappa}u_0(x)=a(t)$ for all $x$ and $t>0$.
The heat equation implies $a'(t)=0$, so $a(t)$ is a constant $a$.
For $\varphi\in C_c^\infty(\R^N)$, the local $L^2$ trace in Lemma~\ref{l:init} implies
$$
 \int_{\R^N}u_0\Delta_\kappa\varphi\,d\mu_\kappa
 =\lim_{t\to0+}\int_{\R^N}e^{t\Delta_\kappa}u_0\Delta_\kappa\varphi\,d\mu_\kappa
 =a\int_{\R^N}\varphi\,d\mu_\kappa.
$$
Thus $\Delta_\kappa u_0=a$ in weighted distributions.
Conversely, if $\Delta_\kappa u_0=c$, testing this identity against
$\chi_R(y)h_t(x,y)$ and letting $R\to\infty$, we obtain
$$
 \Delta_\kappa e^{t\Delta_\kappa}u_0(x)
 =\int_{\R^N}u_0(y)\Delta_{\kappa,y}h_t(x,y)\,d\mu_\kappa(y)
 =c\int_{\R^N}h_t(x,y)\,d\mu_\kappa(y)=c.
$$
Here the cutoff product formula \eqref{e:prod} and
\eqref{e:itail}, also for one derivative of the kernel, make the cutoff errors tend to zero.

\end{remark}

\begin{remark}\label{r:init}
For nonzero initial data, Corollary~\ref{c:en} applies to $u-e^{t\Delta_\kappa}u_0$.
The gradient estimate \eqref{e:en} can fail for $u$ itself.
Indeed, let $f=0$ and choose any nonzero $u_0\in C_c^\infty(\R^N)$.
The initial-value assumptions hold because its Dunkl derivatives are bounded and compactly supported.
There is a ball $B$ on which the mean-square oscillation of $\nabla u_0$ is positive.
Integration of the bound in Lemma~\ref{l:heat}, using doubling, shows that
$$
 \sup_x\int_{\R^N}\|\nabla_xh_s(x,y)\|\,d\mu_\kappa(y)\leq Cs^{-1/2}.
$$
Apply this to the bounded, compactly supported function $\Delta_\kappa u_0$:
\[
 \|\nabla e^{t\Delta_\kappa}u_0-\nabla u_0\|_\infty
 = \left\|\int_0^t\nabla e^{s\Delta_\kappa}\Delta_\kappa u_0\,ds\right\|_\infty
 \leq C\|\Delta_\kappa u_0\|_\infty\int_0^t s^{-1/2}\,ds
 \leq C\sqrt t\,\|\Delta_\kappa u_0\|_\infty\to0.
\]
The oscillation of $\nabla u_0$ on $B$ is positive, so for some $c_B,\varepsilon>0$,
$$
 \frac1{\mu_\kappa(B)}\int_B
 \|\nabla e^{t\Delta_\kappa}u_0-(\nabla e^{t\Delta_\kappa}u_0)_B\|^2\,d\mu_\kappa
 \geq c_B,\qquad 0<t<\varepsilon.
$$
Therefore
$$
 \Xnorm{t^{-1/2}\nabla e^{t\Delta_\kappa}u_0}^2
 \geq c_B\int_0^\varepsilon\frac{dt}{t}=\infty.
$$
\end{remark}

\section{Reflection effects and half-line equations}\label{sec:6}

Throughout this section, $N=1$, $R=\{-1,1\}$ and $\kappa>0$.
Thus $d\mu_\kappa(x)=|x|^{2\kappa}\,dx$ and
\begin{equation}\label{e:rank}
 \Delta_\kappa v(x)
 =v''(x)+\frac{2\kappa}{x}v'(x)-\frac{\kappa}{x^2}\bigl(v(x)-v(-x)\bigr),\qquad x\ne0.
\end{equation}

The term involving $v(-x)$ couples the two half-lines.

\subsection{Response at a reflected point}\label{s:resp}

The slow heat-kernel decay in reflected directions was observed in \cite[Remark~2.4]{ABDH}.

\begin{proposition}\label{p:resp}
Let $a>0$ and $0\leq\varphi\in C_c^\infty((-5a/4,-3a/4))$ with $\varphi(-a)=1$.
Let $u_\kappa$ and $u_{\mathrm{cl}}$ be the Duhamel solutions of the Dunkl and classical heat equations, respectively,
with zero initial value and external force $f(x,t)=\mathbf1_{(0,1)}(t)\varphi(x)$. For $0<t<1$,
\begin{equation}\label{e:resp}
 \left|u_\kappa(a,t)-\frac{\kappa t^2}{2a^2}\right|
 \leq\frac{t^3}{6}\|\Delta_\kappa^2\varphi\|_\infty,
\end{equation}
whereas
$$
 0\leq u_{\mathrm{cl}}(a,t)
 \leq C\|\varphi\|_\infty\,t\exp\left(-\frac{ca^2}{t}\right)
$$
with absolute constants $C,c>0$.
\end{proposition}

\begin{proof}
Since $\varphi$ vanishes near $a$, \eqref{e:rank} implies
$$
 \varphi(a)=\varphi'(a)=\varphi''(a)=0,\qquad
 \Delta_\kappa\varphi(a)=\frac{\kappa}{a^2}.
$$
The Dunkl heat semigroup is a strongly continuous contraction semigroup on $C_0(\R)$,
whose generator agrees with $\Delta_\kappa$ on the Schwartz class; see \cite[Theorem~4.8]{Ro}.
Both $\varphi$ and $\Delta_\kappa\varphi$ are smooth and compactly supported.
Twice integrating the semigroup identity, we obtain
$$
 e^{\tau\Delta_\kappa}\varphi
 =\varphi+\tau\Delta_\kappa\varphi
 +\int_0^\tau(\tau-s)e^{s\Delta_\kappa}\Delta_\kappa^2\varphi\,ds.
$$
Hence, for $0<t<1$,
\[
 u_\kappa(\cdot,t)
 = \int_0^t e^{\tau\Delta_\kappa}\varphi\,d\tau
 = t\varphi+\frac{t^2}{2}\Delta_\kappa\varphi
 +\frac12\int_0^t(t-s)^2e^{s\Delta_\kappa}\Delta_\kappa^2\varphi\,ds.
\]
Evaluating at $x=a$ and using contractivity, we obtain \eqref{e:resp}.

For the classical semigroup, $|a-y|\geq7a/4$ on the support of $\varphi$, so
\[
 0\leq(e^{\tau\Delta}\varphi)(a)
 \leq \frac{\|\varphi\|_\infty}{\sqrt{4\pi\tau}}
 \int_{|a-y|\geq7a/4}\exp\left(-\frac{|a-y|^2}{4\tau}\right)\,dy
 \leq \sqrt2\,\|\varphi\|_\infty\exp\left(-\frac{49a^2}{128\tau}\right).
\]
Integration over $0<\tau<t$ proves the remaining estimate.
The proof is complete.
\end{proof}

In relative form, the Dunkl estimate \eqref{e:resp} reads
$$
 \left|\frac{2a^2u_\kappa(a,t)}{\kappa t^2}-1\right|
 \leq\frac{ta^2}{3\kappa}\|\Delta_\kappa^2\varphi\|_\infty,
 \qquad 0<t<1.
$$
The external force satisfies
$\Xnorm{f}=\|\varphi\|_{\mathrm{BMO}_2(\mu_\kappa)}$, so Theorem~\ref{t:main} yields
$$
 \Xnorm{\partial_tu_\kappa}+\Xnorm{\Delta_\kappa u_\kappa}
 \leq C\|\varphi\|_{\mathrm{BMO}_2(\mu_\kappa)}.
$$

\subsection{A coupled system and an inverse-square potential}\label{s:half}

For $r>0$, write $u_\pm(r,t)=u(\pm r,t)$ and $f_\pm(r,t)=f(\pm r,t)$.
By \eqref{e:rank}, the Dunkl heat equation becomes
\begin{equation}\label{e:cou}
 \partial_tu_\pm-\partial_r^2u_\pm-\frac{2\kappa}{r}\partial_ru_\pm
 +\frac{\kappa}{r^2}(u_\pm-u_\mp)=f_\pm.
\end{equation}
The reflection terms cancel on adding the equations. Subtracting them and writing
$w=u_+-u_-$ and $g=f_+-f_-$, we obtain
\begin{equation}\label{e:odd}
 \partial_tw-\partial_r^2w-\frac{2\kappa}{r}\partial_rw
 +\frac{2\kappa}{r^2}w=g,\qquad r,t>0,\qquad w(\cdot,0)=0.
\end{equation}

We extend a function $q$ on $(0,\infty)\times\R_+$ to the whole line in two ways:
$$
 q^{\mathrm{odd}}(x,t)=\operatorname{sgn}(x)q(|x|,t),\qquad
 q^{\mathrm{even}}(x,t)=q(|x|,t),\qquad x\ne0.
$$
In terms of intervals on the half-line, \eqref{e:bmo} yields the comparison
\begin{align*}
 \Xnorm{q^{\mathrm{odd}}}^2\simeq&
 \sup_{I\subset(0,\infty)}\frac1{\mu_\kappa(I)}
 \int_0^\infty\int_I|q(r,t)-(q(\cdot,t))_I|^2\,d\mu_\kappa(r)\,dt\\
 &+\sup_{\rho>0}\frac1{\mu_\kappa((0,\rho))}
 \int_0^\infty\int_0^\rho|q(r,t)|^2\,d\mu_\kappa(r)\,dt.
\end{align*}
Here $I$ ranges over bounded intervals and $d\mu_\kappa(r)=r^{2\kappa}\,dr$.
Restricting the supremum in \eqref{e:bmo} to intervals in $(0,\infty)$ bounds the first term by $\Xnorm{q^{\mathrm{odd}}}^2$.
Taking the interval $(-\rho,\rho)$ bounds the second, since the odd extension has mean zero there.
Conversely, by reflection the first term controls intervals in either half-line.
For $J=(-a,b)$ with $a,b>0$ and $\rho=\max\{a,b\}$,
\begin{align*}
 \frac1{\mu_\kappa(J)}\int_0^\infty\int_J
 |q^{\mathrm{odd}}-(q^{\mathrm{odd}}(\cdot,t))_J|^2\,d\mu_\kappa\,dt
 & \leq\frac1{\mu_\kappa(J)}\int_0^\infty\int_J|q^{\mathrm{odd}}|^2\,d\mu_\kappa\,dt\\
 & \leq\frac2{\mu_\kappa((0,\rho))}\int_0^\infty\int_0^\rho|q|^2\,d\mu_\kappa\,dt.
\end{align*}

For the system \eqref{e:cou}, reflection invariance of $\mu_\kappa$ implies
$$
 g^{\mathrm{odd}}(x,t)=f(x,t)-f(-x,t),\qquad
 \Xnorm{g^{\mathrm{odd}}}\leq2\Xnorm{f}.
$$

By a weak solution of \eqref{e:odd}, we mean a function $w$ whose odd extension
satisfies \eqref{e:weak} on $\R$, with external force $g^{\mathrm{odd}}$ and initial value zero.
The half-line operator is obtained by restricting $-\Delta_\kappa$ to odd functions.
Its heat kernel, relative to $y^{2\kappa}\,dy$ on $(0,\infty)$, is
$h_t(r,y)-h_t(r,-y)$.

\begin{corollary}\label{c:odd}
Let $g$ be a half-line external force with $\Xnorm{g^{\mathrm{odd}}}<\infty$. The function
$$
 w(r,t)=\int_0^t\int_0^\infty
 \bigl[h_{t-s}(r,y)-h_{t-s}(r,-y)\bigr]g(y,s)y^{2\kappa}\,dy\,ds
$$
is the unique weak solution of \eqref{e:odd} whose odd extension satisfies \eqref{e:class}. Moreover,
$$
 \Xnorm{(\partial_tw)^{\mathrm{odd}}}
 +\biggl\|\biggl(\partial_r^2w+\frac{2\kappa}{r}\partial_rw-\frac{2\kappa}{r^2}w\biggr)^{\mathrm{odd}}\biggr\|_{\mathbb L^2(\R_+;\mathrm{BMO}(\mu_\kappa))}
 \leq C\Xnorm{g^{\mathrm{odd}}},
$$
and
$$
 \Xnorm{t^{-1/2}(\partial_rw)^{\mathrm{even}}}
 +\Xnorm{t^{-1/2}(w/r)^{\mathrm{even}}}
 \leq C\Xnorm{g^{\mathrm{odd}}}.
$$
\end{corollary}

\begin{proof}
The average $(g^{\mathrm{odd}}(\cdot,t))_{(-\rho,\rho)}$ vanishes by oddness. Hence
$$
 \int_0^T\int_{-\rho}^{\rho}|g^{\mathrm{odd}}(x,t)|^2\,d\mu_\kappa(x)\,dt
 \leq\mu_\kappa((-\rho,\rho))\Xnorm{g^{\mathrm{odd}}}^2,\qquad \rho,T>0.
$$
This verifies the local integrability required in Theorem~\ref{t:main}. Its solution
$$
 u(x,t)=\int_0^t e^{(t-s)\Delta_\kappa}g^{\mathrm{odd}}(\cdot,s)(x)\,ds
$$
is odd because reflection commutes with the heat semigroup.
Changing $y$ to $-y$ on the negative half-line identifies $u$ with $w^{\mathrm{odd}}$.
On each half-line,
$$
 \Delta_\kappa u
 =\left(\partial_r^2w+\frac{2\kappa}{r}\partial_rw-\frac{2\kappa}{r^2}w\right)^{\mathrm{odd}}
$$
in distributions. The equation and maximal regularity estimate therefore follow from Theorem~\ref{t:main},
as does uniqueness among solutions with odd extension in \eqref{e:class}.
Differentiating the odd extension, we find
$$
 \partial_xu=(\partial_rw)^{\mathrm{even}}.
$$
The reflection term satisfies
$$
 \frac{u(x,t)-u(-x,t)}{x}=2(w/r)^{\mathrm{even}}(x,t).
$$
Apply \eqref{e:en} to the first identity and \eqref{e:ren} to the second.
The proof is complete.
\end{proof}

\appendix
\section{Direct estimates for short and long times}\label{sec:a}

Fix $B=B(c,r)$. We prove Proposition~\ref{p:mr} and the upper bounds in Theorem~\ref{t:sq}
by splitting the integrals in $\tau$ at $r^2$. For maximal regularity, write
$$
 f_0=\bigl(f-(f(\cdot,t))_{4B}\bigr)\mathbf1_{4B},\qquad
 f_\infty=\bigl(f-(f(\cdot,t))_{4B}\bigr)\mathbf1_{(4B)^c},
$$
and set
$$
 V_\tau(x)=(\Delta_\kappa e^{\tau\Delta_\kappa}f_\infty)(x,\cdot-\tau),
$$
where $f_\infty$ is extended by zero for negative times. By \eqref{e:l2} and doubling,
$$
 \frac1{\mu_\kappa(B)}\int_B
 \|\mathcal M_\kappa^{a,b}f_0(x,\cdot)\|_{L^2(\R_+,dt)}^2\,d\mu_\kappa(x)
 \leq\frac{4\|f_0\|_{L^2(\mu_\kappa\times\R_+)}^2}{\mu_\kappa(B)}\leq C\Xnorm{f}^2.
$$
For the other part, take $x\in B$ and $y\notin4B$.
 Since $\|x-y\|\geq3r$, for $0<\tau<r^2$ we have
$$
 \left(1+\frac{\|x-y\|}{\sqrt\tau}\right)^{-2}
 \leq C\frac{\tau}{r^2}
       \left(1+\frac{\|x-y\|}{r}\right)^{-2}.
$$
After multiplication by $\tau^{-1}$ from \eqref{e:heat}, the right-hand side is $Cr^{-2}$ times
the Euclidean weight in Lemma~\ref{l:mean} with $\eta=2$, since $\max\{r,\sqrt\tau\}=r$.
When $\tau\geq r^2$, subtract $V_\tau(c)$, which does not affect spatial oscillation,
and integrate the spatial derivative along $[c,x]$. This gives an extra factor $r/\sqrt\tau$.
Minkowski's inequality and contractivity of time translation, followed by Lemmas~\ref{l:heat}
and \ref{l:mean}, therefore imply
\begin{align*}
 \left\|\int_{(a,b)\cap(0,r^2)}V_\tau(x)\,d\tau\right\|_{L^2(\R_+)}
 & \leq C\Xnorm{f}\int_0^{r^2}\frac1\tau\frac\tau{r^2}
      \left(1+\log\frac r{\sqrt\tau}\right)\,d\tau\\
 &\leq C\Xnorm{f},
\end{align*}
and
\begin{align*}
 \left\|\int_{(a,b)\cap(r^2,\infty)}[V_\tau(x)-V_\tau(c)]\,d\tau\right\|_{L^2(\R_+)}
 & \leq C\Xnorm{f}\int_{r^2}^\infty\frac1\tau\frac r{\sqrt\tau}
      \left(1+\log\frac{\sqrt\tau}r\right)\,d\tau\\
& \leq C\Xnorm{f}.
\end{align*}
For fixed $a,b$, the integral of $V_\tau(c)$ over $(a,b)\cap(r^2,\infty)$
belongs to $L^2(\R_+,dt)$ by Lemma~\ref{l:mean}.
Subtracting it and then centering by the spatial average proves Proposition~\ref{p:mr}.

For the square functions, write $F-F_{4B}=F_0+F_\infty$, where
$$
 F_0=(F-F_{4B})\mathbf1_{4B},
 \qquad
 F_\infty=(F-F_{4B})\mathbf1_{(4B)^c}.
$$
Fix one of the two operator families
$$
 F\longmapsto\nabla e^{\tau\Delta_\kappa}F,
 \qquad
 F\longmapsto
 \frac{e^{\tau\Delta_\kappa}F
 -(e^{\tau\Delta_\kappa}F)\circ\sigma_\alpha}
 {\langle\alpha,\cdot\rangle}
 \quad(\alpha\in R_+,\ \kappa(\alpha)>0),
$$
and denote it by $\mathcal T_\tau$, $\tau>0$.
Both operators vanish on spatial constants, so
$$
 \mathcal T_\tau F
 =\mathcal T_\tau F_0+\mathcal T_\tau F_\infty.
$$
By \eqref{e:el2} and doubling,
\begin{equation}\label{e:aloc}
 \frac1{\mu_\kappa(B)}
 \int_B\int_0^\infty
 \|\mathcal T_\tau F_0(x)\|^2\,d\tau\,d\mu_\kappa(x)
 \leq
 \frac{C\|F_0\|_{L^2(\mu_\kappa;\mathcal H)}^2}
 {\mu_\kappa(B)}
 \leq C\Bnorm{F}^2.
\end{equation}

For $x\in B$, $y\notin4B$, and $0<\tau<r^2$,
$$
 \left(1+\frac{\|x-y\|}{\sqrt\tau}\right)^{-1}
 \leq
 C\frac{\sqrt\tau}{r}
 \left(1+\frac{\|x-y\|}{r}\right)^{-1}.
$$
Thus Lemmas~\ref{l:heat}, \ref{l:ker}, and \ref{l:mean},
with $\eta=1$ in the last lemma, imply
$$
 \|\mathcal T_\tau F_\infty(x)\|
 \leq
 \frac Cr\left(1+\log\frac r{\sqrt\tau}\right)\Bnorm{F},
 \qquad x\in B,\quad 0<\tau<r^2.
$$
For $\tau\geq r^2$, the kernel and spatial derivative bounds imply
$$
 \|\mathcal T_\tau F_\infty(c)\|
 +\sqrt\tau\sup_{z\in B}
 \|\nabla\mathcal T_\tau F_\infty(z)\|
 \leq
 \frac C{\sqrt\tau}
 \left(1+\log\frac{\sqrt\tau}r\right)\Bnorm{F}.
$$
Integrating the spatial derivative along the segment from $c$ to $x$,
we obtain
\begin{align*}
 \|\mathcal T_\tau F_\infty(x)
       -\mathcal T_\tau F_\infty(c)\|
 &\leq
 \|x-c\|\int_0^1
 \|\nabla\mathcal T_\tau F_\infty(c+\theta(x-c))\|\,d\theta\\
 &\leq
 \frac{Cr}{\tau}
 \left(1+\log\frac{\sqrt\tau}r\right)\Bnorm{F},
 \qquad x\in B,\quad \tau\geq r^2.
\end{align*}
Consequently, for every $x\in B$,
\begin{align}\label{e:afar}
 &\int_0^{r^2}\|\mathcal T_\tau F_\infty(x)\|^2\,d\tau
 +\int_{r^2}^\infty
 \|\mathcal T_\tau F_\infty(x)
       -\mathcal T_\tau F_\infty(c)\|^2\,d\tau
 \nonumber\\
 \leq\,&
 C\Bnorm{F}^2
 \left\{
 \int_0^{r^2}\frac1{r^2}
 \left(1+\log\frac r{\sqrt\tau}\right)^2\,d\tau
 +\int_{r^2}^\infty\frac{r^2}{\tau^2}
 \left(1+\log\frac{\sqrt\tau}r\right)^2\,d\tau
 \right\}
 \nonumber\\
 =  \,&
 C\Bnorm{F}^2
 \left\{
 \int_0^1\left(1-\frac12\log v\right)^2\,dv
 +\int_1^\infty
 \frac{\left(1+\frac12\log v\right)^2}{v^2}\,dv
 \right\}
 \leq C\Bnorm{F}^2.
\end{align}

For $0<a<b<\infty$, the estimate at $c$ for $\tau\geq r^2$ implies that
$\mathbf1_{[r^2,\infty)}(\tau)\mathcal T_\tau F_\infty(c)$ is square-integrable on $(a,b)$.
By the minimizing property of the spatial average, \eqref{e:aloc} and \eqref{e:afar} imply
\begin{align*}
 &\frac1{\mu_\kappa(B)}
 \int_B\int_a^b
 \|\mathcal T_\tau F(x)-(\mathcal T_\tau F)_B\|^2
 \,d\tau\,d\mu_\kappa(x)\\
 \leq\,&
 \frac1{\mu_\kappa(B)}
 \int_B\int_a^b
 \bigl\|
 \mathcal T_\tau F_0(x)+\mathcal T_\tau F_\infty(x)
 -\mathbf1_{[r^2,\infty)}(\tau)\mathcal T_\tau F_\infty(c)
 \bigr\|^2\,d\tau\,d\mu_\kappa(x)\\
 \leq\,& C\Bnorm{F}^2.
\end{align*}
The bound is uniform in $a,b$ and $B$. By monotone convergence it holds on $(0,\infty)$.
Taking the supremum over $B$ proves \eqref{e:sq} and \eqref{e:rsq}.
Since
$$
 \nabla_\kappa e^{\tau\Delta_\kappa}F
 =
 \nabla e^{\tau\Delta_\kappa}F
 +\sum_{\substack{\alpha\in R_+\\\kappa(\alpha)>0}}
 \kappa(\alpha)\alpha
 \frac{e^{\tau\Delta_\kappa}F
 -(e^{\tau\Delta_\kappa}F)\circ\sigma_\alpha}
 {\langle\alpha,\cdot\rangle},
$$
the triangle inequality also proves the upper bound in \eqref{e:dsq}.

By \eqref{e:aloc} and the integral over $(0,r^2)$ in \eqref{e:afar},
\begin{align*}
 \frac1{\mu_\kappa(B)}
 \int_B\int_0^{r^2}
 \|\mathcal T_\tau F(x)\|^2\,d\tau\,d\mu_\kappa(x)
 &\leq\frac2{\mu_\kappa(B)}
 \int_B\int_0^{r^2}
 \bigl(\|\mathcal T_\tau F_0(x)\|^2+\|\mathcal T_\tau F_\infty(x)\|^2\bigr)\,d\tau\,d\mu_\kappa(x)\\
 &\leq C\Bnorm{F}^2.
\end{align*}
Combining the gradient and reflection terms in
$\nabla_\kappa e^{\tau\Delta_\kappa}F$ and taking the supremum over $B$,
we conclude that
\begin{equation}\label{e:short}
 \sup_B\frac1{\mu_\kappa(B)}
 \int_0^{r_B^2}\int_B
 \|\nabla_\kappa e^{\tau\Delta_\kappa}F(x)\|_{\mathcal H^N}^2
 \,d\mu_\kappa(x)\,d\tau
 \leq C\Bnorm{F}^2.
\end{equation}

\subsection*{Acknowledgements}\label{s:ack}
The authors would like to thank Prof. Ji Li  for helpful discussions and comments.  During the preparation of the manuscript, OpenAI’s ChatGPT was used to assist
with checking proofs and language editing. All mathematical arguments and
references were independently verified by the authors, who take full
responsibility for the content of the paper.

\end{document}